\documentclass[12pt,a4paper]{article}
\usepackage{amsmath}
\usepackage{amssymb}
\advance\topmargin by -2.2cm                    
\newtheorem{theo}{\indent Theorem}[section]
\newtheorem{prop}[theo]{\indent Proposition}

\newtheorem{rem}[theo]{\indent Remark}
\newtheorem{lem}[theo]{\indent Lemma}

\newtheorem{cor}[theo]{\indent Corollary}

\newtheorem{ex}[theo]{\indent Example}
\newtheorem{ass}[theo]{\indent Assumption}

\newenvironment{proof}{\noindent {\bf Proof }}
{\hfill $\bullet$ \vspace{0.25cm}}

\def\R{\mathbb{R}}
\def\1{\mathbf{1}}
\def \E{\mathbb{E}}
\def \P{\mathbb{P}}
\def \N{\mathbb{N}}
\def \L{\mathbb{L}}

\def\F{\cal{F}}
\def \ge{\varepsilon}

\begin{document}
\newlength{\breite}
\settowidth{\breite}{Web: http://www.mathematik.uni-mainz.de/Stochastik/Loecherbache}

\title{Polynomial bounds in the Ergodic Theorem for 
one-dimensional diffusions and integrability of hitting times}

\author{Eva 
{\sc L\"ocherbach}\footnote{Centre de Math\'ematiques,
Facult\'e de Sciences et Technologie,
Universit\'e Paris XII, 61 avenue du G\'en\'eral de Gaulle, 94010 Cr\'eteil 
Cedex, France. E-mail: {\tt 
locherbach@univ-paris12.fr}} , Dasha {\sc Loukianova}\footnote{D\'epartement de Math\'ematiques, Universit\'e d'Evry-Val d'Essonne, Bd Fran\c{c}ois Mitterrand, 91025 Evry Cedex, France. E-mail: {\tt 
 dasha.loukianova@univ-evry.fr
} },
Oleg {\sc Loukianov}\footnote{D\'epartement Informatique, IUT S\'enart-Fontainebleau, Universit\'e Paris XII, route Hurtault, 77300 Fontainebleau, France. E-mail: {\tt 
oleg.loukianov@univ-paris12.fr}}}

\maketitle

\def\abstractname{Abstract}
\begin{abstract}
Let $X$ be a one dimensional positive recurrent diffusion with initial distribution $\nu$ and invariant probability $\mu$. Suppose that for some $p> 1$, $\exists a\in\R$ such that $\forall x\in\R,\ \E_x T_a^p<\infty$ and $\E_\nu T_a^{p/2}<\infty$, where $T_a$ is the hitting time of $a$.  For such a diffusion, we derive non asymptotic deviation bounds of the form 
$$\P_{\nu}\left (\left|\frac1t\int_0^tf(X_s)ds-\mu(f)\right|\geq\ge\right)\leq K(p)\frac1{t^{p/2}}\frac 1{\ge^p}A(f)^p.$$
Here $f$  bounded or bounded and compactly supported and $A(f)=\|f\|_{\infty}$ when $f$ is bounded and $A(f)=\mu(|f|)$ when $f$ is bounded and compactly supported.

We also give, under some conditions on the coefficients of $X$, a polynomial control of $\E_xT_a^p$ from above and below. This control is  based on a generalized Kac's formula (see theorem~\ref{thm:mainKac}) for the moments $\E_x f(T_a)$ of a differentiable function $f$.\\
\begin{center}
{\bf R\'esum\'e}
\end{center}
Consid\'erons une diffusion r\'ecurrente positive avec loi initiale $\nu $ et probabilit\'e invariante $\mu .$ Pour tout $a \in \R,$ soit $T_a$ le temps d'atteinte du point $a.$ Supposons qu'il existe $p > 1 $ et un point $ a \in \R $ tels que pour tout $x \in \R,$ $ \E_x T_a^p < \infty $ et $\E_\nu T_a^{p/2}<\infty.$ Alors nous obtenons l' in\'egalit\'e de d\'eviation non-asymptotique suivante : 
$$\P_{\nu}\left (\left|\frac1t\int_0^tf(X_s)ds-\mu(f)\right|\geq\ge\right)\leq K(p)\frac1{t^{p/2}}\frac 1{\ge^p}A(f)^p,$$
o\`u $f$ est une fonction born\'ee ou une fonction born\'ee \`a support compact. Ici, $A(f)=\|f\|_{\infty}$ dans le cas d'une fonction born\'ee et $A(f)=\mu(|f|) $ dans le cas d'une fonction born\'ee \`a support compact. 

De plus, sous certaines conditions sur les coefficients de la diffusion, nous obtenons une minoration et majoration, polynomiale en $x,$ de $\E_xT_a^p .$ Ce r\'esultat est bas\'e sur une formule de Kac g\'en\'eralis\'ee (voir th\'eor\'eme \ref{thm:mainKac}) pour les moments $\E_x f(T_a)$ o\`u $f$ est une fonction d\'erivable.        
\end{abstract}


{\it MSC 2000}  : 60 F 99, 60 J 55, 60 J 60

\section{Introduction}

We consider the solution of the  one-dimensional stochastic differential equation
\[dX_t=\beta (X_t)\,dt+\sigma(X_t)\,dW_t , \]
with arbitrary  initial data. Suppose that $X$ is positive recurrent, and denote by $\mu$ its invariant probability. From the Ergodic Theorem in this case we know that for all $x\in\R $, $f\in\L^1(\mu)$ and $\ge>0$ 
 \begin{equation}\label{ergthm}
  \P_x\left (\left|\frac1t\int_0^tf(X_s)ds-\mu(f)\right|\geq\ge\right)\to 0
  \end{equation}
  as $t$ goes to $+\infty.$ 
 The purpose of this paper is to obtain a non-asymptotic upper bound for the probability in \eqref{ergthm}. 
%Moreover, we try to impose rather week condition on $X$, hence we restrict ourself to the study
%of the case when the speed of convergence in \eqref{ergthm} is polynomial.
Such a bound  is of major importance for many applications: various non asymptotic problems for statistics of diffusions (see \cite{comteetal}, \cite {GaP2}, \cite{drift}), concentration for particular approximations of granular media equations (see \cite{CaGuiMal}), and many other examples. Mainly, such a bound is useful any time when we wish to substitute a random quantity $\frac1t\int_0^tf(X_s)ds$ by a deterministic $\mu(f)$ except on some set of ``small'' probability. ``Small'' usually means ``exponentially small'', and this case  has already been discussed in the literature, (see the references below). Other possible rates seem not to be studied so far, but actually it turns out that in concrete problems it is often sufficient to consider slower rates. On the other hand, considering slower rates generally permits to lighten the assumptions on the model.  
 
In this paper we study the case when the rate of convergence in \eqref{ergthm} is polynomial. We use the regeneration method, which appeals to the following natural condition: the integrability of regeneration times. For bounded or bounded and compactly supported functions $f$, and $X$ such that for some $p> 1$  the $p$-th moment of the regeneration time exists, we show the following deviation inequality:  for all $0<\ge<A(f),$
\begin{equation}\label{erginequality}
\P_{\nu}\left(\left|\frac1t\int_0^tf(X_s)ds-\mu(f)\right|>\ge \right)\leq K( p, x)\ge^{-p}t^{-\alpha/2}A(f)^p.
\end{equation}
Here
 $A(f)=\|f\|_{\infty}$ when $f$ is bounded and $A(f)=\mu(|f|)$ when $f$ is bounded and of compact support, $\alpha=p$ if $p\geq 2$, and $\alpha=p-1$ if $1<p\leq 2$. The constant $K$ is a positive constant, which does not depend on $f$, $t$, $\ge,$ see  Theorem \ref{mainbounded}, Theorem \ref{mainint} for the precise statement.

Since the one-dimensional case is very explicit, the moments of regeneration times are closely related to hitting time moments. In the last section we formulate conditions on the existence of hitting time moments and give (corollary \ref{expbounded}, \ref{explintegrable}) some sufficient conditions for \eqref{erginequality}  in terms of the coefficients of the diffusion.
   %One advantage of the regeneration method is that in contrary to \cite{lezaud} or \cite{CaGui}, where only $\mu$-absolutely continuous starting measures can be used, our inequality is valuable for any initial distribution.

Let us give a short overview of the history of the problem.
In the context of i.i.d.~variables the question of the rate of convergence in \eqref{ergthm} turns out to be the question of the rate of convergence in the Law of Large Numbers. This rate is exponential whenever the variables have exponential moments. There is a large literature on this subject;  let us cite very early results by  Bernstein, \cite {Pe},
 Bennet, \cite{Be}, Hoeffding, \cite{Hoef}, the book by Petrov, \cite {Pe}, a more recent article by Pinelis, \cite {Pin}, and the references therein. 
 
For Markov chains, Cl\'emen\c{c}on, \cite{Cl}, deduces an exponential bound for the probability in \eqref {ergthm} using the regeneration method. He works with geometrically regular Markov chains, which means exponential integrability of some hitting times, in the stationary regime and with bounded functions $f$ (see also Bertail-Cl\'emen\c{c}on \cite{BeCl}). 
  Following a completely different approach, Adamczak, \cite {Adam}, derives concentration inequalities for empirical processes of Markov chains. As a particular case he deduces an exponential bound for \eqref{ergthm}, when $f$ is bounded and $\mu(f)=0.$ He also works under the assumption of exponential integrability of some regeneration time.  As far as we know, in the context of Markov chains the polynomial rate of convergence in the Ergodic theorem has not been considered. However, it has been studied for many other ergodic phenomena, see for example Tuominen and Tweedy, \cite{TuoTw}, Jarner and Roberts \cite{JaRo}, Chazottes and Redig, \cite{chazottes} and references therein. Moreover, it is a well-known observation that there is a natural connection between the speed of convergence to equilibrium and the integrability of some stopping (typically regeneration or coupling) times, see Chazottes and Redig, \cite {chazottes}, Meyn and Tweedie, \cite{mt}, chapter III, 15 , 
Douc-Fort-Moulines-Soulier, \cite{DFMS}.   

For continuous time Markov processes, as we already mentioned, the non-asymptotic bound in the Ergodic theorem  was obtained  by Lezaud, \cite{lezaud}, and Cattiaux and Guillin, \cite{CaGui}. The approach in \cite{CaGui}  relies on the use of functional inequalities for the invariant probability $\mu$ like the Poincar\'e inequality. In this way the authors obtain an asymptotically sharp exponential bound, in the spirit of the large deviation principle (see also \cite{wu}), for a process starting from the invariant measure $\mu$ or from an initial law being absolute continuous with respect to $\mu$. Another approach is followed in \cite {lezaud} where perturbation of operator theory is used. All these authors work under the assumption of a spectral gap and obtain an exponential bound for \eqref{ergthm}. Concerning the spectral gap, recently Loukianov, Loukianova and Song, \cite{LLS}, proved that 
this condition is equivalent to the existence of exponential moments of hitting times for one dimensional diffusions. 
%Hence we observe for another time the relation between the rate of convergence to equilibrium and the integrability of regeneration times.  
Note also that more general exponential bounds are obtained in Guillin, L\'eonard, Wu and Yao, \cite{glwy}, in relation with transportation of measure inequalities. 
For one dimensional ergodic diffusion processes, Galtchouk and Pergamenshchikov, \cite{GaP}, obtain \eqref{ergthm} uniformly with respect to the initial condition and to some other parameter. Their bound is exponential, too. They work under the assumption of constant diffusion coefficient and a drift bounded from above and below by linear functions.
Finally let us also mention the paper of Kontoyiannis and Meyn, \cite{km1}, where an exponential bound for the integral version of \eqref{ergthm} is obtained. This work concerns multiplicatively (and geometrically) regular  Markov processes, see also \cite{km2} for its discrete counterpart. 

In the continuous time, to the best of our knowledge, the polynomial case of the rate of convergence in the Ergodic theorem  \eqref{ergthm} has not yet been considered. However there are a lot of results on polynomial rates for other phenomena of convergence to equilibrium. The most studied are the rate of decrease of mixing coefficients and the rate of decrease of the total variation distance between the law of $X_t$ and $\mu  . $ When the last rate is exponential (resp.~sub-exponential or polynomial), the model is usually called exponentially (resp.~sub-exponentially or polynomially) ergodic. In this field of research,
Fort and Roberts, \cite{FR}, study the sub-exponential ergodicity for a strong Markov process
and obtain as an application of their results the polynomial ergodicity for  multi-dimensional diffusions.  Veretennikov, \cite{Veret}, studies both mixing coefficients and total variation distance between the law of $X_t$ and $\mu $ and gives sufficient conditions for their polynomial decrease in the framework of  multi-dimensional diffusions. The conditions in \cite{FR} and \cite{Veret} are formulated in terms of the coefficients of the diffusion, but both papers involve the existence of polynomial  moments for some regeneration times: modulated for \cite{FR} and coupling times for \cite{Veret}. Finally, Douc, Fort and Guillin, \cite{DFG}, study sub-geometric ergodicity of a strong Markov process and provide a criterion that yields a precise control of a sub-geometric moment of the return-time to a test-set (modulated moment). 
Hence the  relation between the integrability of regeneration times and  different types of ergodicity in the sense of total variation distance between the law of $X_t$ and $\mu$ seems to be quite well understood. Regarding the very huge literature on this subject, let us also cite Roberts and Tweedie, \cite {RoTw}, Down, Meyn and Tweedie, \cite{DMT}, Douc, Guillin and Moulines, \cite {DouGuiMou}, Pardoux and Veretennikov, \cite{Pardoux-Veret}, Veretennikov and Klokov, \cite{Veret-Klok}, and the references therein.
 
 In this paper, we establish a very explicit relation between integrability of hitting times and speed of convergence in the Ergodic theorem \eqref{ergthm}. Hence a large part of the paper is devoted to the study of hitting time's moments. 
In Theorem \ref{indepofpoint} we explain  that  $\E_xT_y^p$ is finite or infinite simultaneously for all couples $x<y$ or $x>y$. The proof of this result is based on a generalized  version  of Kac's moment formula (Theorem \ref{thm:mainKac}), interesting in its own. Recall that the original Kac's formula given in \cite{FitsKac}  relates the moment of order $p$ of hitting times $T_y $ (or more generally of a stopped additive functional) to the previous moment of order $p-1,$ for any $p \in [1, + \infty[.$  Our version (Theorem \ref{thm:mainKac}) relates the moment $\E_xf(T_y)$ to the moment of $\E_xf'(T_y)$.    

In order to be able to work with an initial distribution $\nu$ and to check $\E_{\nu} T_y^p < \infty ,$ we give in Theorem \ref{theo:control} upper and lower polynomial bounds  for $\E_{x}T_y^p$ under assumptions in the spirit of those given by  Veretennikov, \cite {Veret},  and by Balaji and Ramasubramanian, see \cite {bhatt}. The constants in our bounds are sharp.
A comparative analysis of our conditions with those of \cite {Veret} and \cite{bhatt} is contained in the last section $5$.
 
The paper is organized as follows. Section $2$  collects auxiliary probabilistic results, needed for the proof of the deviation theorems. The Deviations theorems are stated and proved in section $3$. They hold true under the assumption that $\E_x T_y^p < \infty $ for all $x,y .$ Consequently,
sections 4 and 5 are devoted to the study of polynomial integrability of hitting times:
section 4 contains generalized Kac's formula and theoretical conditions for $\E_x T_y^p < \infty $ for all $x,y .$ A precise polynomial bounds for $\E_x T_y^p ,$ under conditions on the coefficient of $X$, as well as  some sufficient conditions for \eqref{erginequality}  in terms of the coefficients of the diffusion are given in the last section. 
 
\vskip0.5cm \noindent 
{\bf Acknowledgments.} The authors  thank Valentine Genon-Catalot and Shiqi Song for very useful discussions concerning hitting times.   

Eva L\"ocherbach has been partially supported by an ANR projet : Ce travail a b\'en\'efici\'e d'une aide de l'Agence Nationale de la Recherche
portant la r\'ef\'erence ANR-08-BLAN-0220-01.
\section{Notation, basic assumptions and auxiliary results}
Let  $X_t$ be a one-dimensional
diffusion process given by
\begin{equation}\label{diffusion}
dX_t=\beta (X_t)\,dt+\sigma(X_t)\,dW_t .
\end{equation}

We impose the following condition on the coefficients of \eqref{diffusion}.

\begin{ass}\label{model}
\begin{enumerate}
\item For all $x, $ $ \sigma^2 (x)>0  .$
\item $\beta $ and $\sigma$ are locally Lipschitz, and $|\sigma(x)|+|\beta(x)|\leq C(1+|x|)$, for some $C>0.$\end{enumerate}
\end{ass}
This assumption ensures the existence of a unique strong non-exploding solution of \eqref{diffusion} (see for example \cite {BS}, Chapter III. 4.17).

 Let us recall some basic facts about one-dimensional diffusions. Denote
$$ s(x)=\exp\left(-2\int_0^x\frac{\beta(u)}{\sigma^2(u)}du\right), \quad  m(x)=\frac{2}{\sigma^2(x)s(x)} ,$$ 
and recall that the scale function is given by
 $$S(x)=\int_0^xs(t)dt \quad\mbox{ for}\quad x \geq 0, \quad  S(x) =  - \int_{x}^0 s(t) dt \quad\mbox{ for}\quad x < 0 .$$
The diffusion $X$ is said to be recurrent if for all $ x\in\R, y\in\R ,$ $ \P_x(T_y<\infty)=1$.
A necessary and sufficient condition of recurrence is
  \begin{equation}\label{rec}
 \lim_{x \to + \infty } S(x) = + \infty \quad \mbox{ and}\quad \lim_{x \to - \infty} S(x) = - \infty
 \end{equation}
  (see \cite{Kh}, example 2 in section 3.8, or \cite {RY},  Ch.VII  ex.$3.21$). 
A recurrent diffusion is called positively recurrent if $
\E_x(T_y)<\infty$ for all $x,y \in \R.$ This condition is equivalent to
$$ M := \int_{-\infty}^\infty m(x) dx < +\infty $$
(see \cite {BS} Chapter II.1.12.).
In the case of positive recurrence, the unique invariant probability measure of the process is given by 
\begin{equation}\label{eq:mu}
 \mu(dx)=\frac 1M m(x) dx .
 \end{equation}
For the remainder of the article, except Proposition \ref{iid}, we suppose 
 \begin{ass}\label{posrec}
  $X$ is positively recurrent.
\end{ass} 

%\section{Auxiliary results.}  
In the sequel we use the regeneration method for one-dimensional diffusions. 
One possible way to introduce the regeneration times is the following: Fix two
points $a<b,$ $ a, b \in \R .$
Define a sequence of stopping times $(S_n)_n$, $(R_n)_n$ as follows :
 $S_0=0$, $R_0=0,$
$$ S_1 := \inf
 \{t\geq 0 : X_t = b \},\quad  R_1 := \inf \{ t \geq S_1 : X_t = a \} ,$$
and for $n \geq 1,$ 
$$S_{n+1} := \inf \{ t > R_n : X_t = b \}, \quad   R_{n+1} := \inf \{ t \geq S_{n+1}  : X_t = a \} . $$
The sequence $(R_n)_n$ ``cuts" the process into i.i.d.~blocs in the following sense: 
If $f:\R\to\R$ is measurable and bounded and if we put
\begin{equation*}%\label{xi}
\xi_n=\int_{R_n}^{R_{n+1}}f(X_s)ds, \quad n\geq 0 ,
\end{equation*}
then we have the following proposition. 

\begin{prop}\label{iid}
Suppose that assumption \ref{model} and condition \eqref{rec} hold. For any initial distribution $\nu$, the sequence $(\xi_n)_{n\geq 1}$ is an i.i.d. sequence under $\P_{\nu}$. For all $n\geq 1,$ the law of $\xi_n$ under $\P_{\nu}$ is equal to the law of $\xi_0$ under $\P_a.$
 \end{prop} 
This last proposition is well known and easy to show using
 the strong Markov property. Note that in particular the sequence $(R_{k+1}-R_{k}), k=1,2,\ldots$ is an i.i.d.~sequence with common distribution equal to the law of $R_1$ under $\P_a$.
%Now suppose that  $f$ has  compact support and is bounded.
Denote
$$ C(f) := \sup_x \E_x \int_0^{R_1} |f|(X_s) ds .$$

\begin{prop}\label{fspe}  
Grant assumptions \ref{model} and \ref{posrec}. If $f$ is measurable bounded with compact support, then
$ C(f)  < \infty .$
\end {prop}

%\begin{rem} In an abstract setting of Markov processes, functions satisfying $C(f)<\infty$ are called special. 
%\end{rem}

\begin{proof}
Denote by $K$ the support of $f$ and let $\tau=\inf\{t\geq 0: \; X_t\in K\}.$
Let $M>0$ be such that $|f|\leq M$. Then,
\begin{multline*}
C(f)=\E_x \left( \E_x\left(\int_{\tau \wedge R_1}^{R_1}|f(X_s)|ds|{\F}_{\tau}\right)\right) \le \E_x\left( \E_{X_\tau}\int_0^{R_1}|f(X_s)|ds\right) \\
\le \sup_{x\in K}\E_x\int_0^{R_1}|f(X_s)|ds\le M\sup_{x\in K}\E_x{R_1} .
\end{multline*}
Since $X$ is positive recurrent, we can use Theorem \ref{indepofpoint} below for $n=1$. This theorem implies that  $x\mapsto \E_x{R_1}$ is continuous, and thus $\sup_{x\in K}\E_x{R_1}<\infty.$
 \end{proof}

Note that the last proposition is true in a much more general case. Actually it is true for any recurrent strong-Feller diffusion  with  state space $\R^n$, see Remark 5.28, 4) of \cite{holo}.
  
The following proposition extends the uniform in $x$ integrability  property of the first  life cycle  and  will play an important role in the sequel.
 
\begin{prop}\label{prop:momentp}
Grant assumptions \ref{model} and \ref{posrec}. Let $f$ be a bounded measurable function with compact support. 
Then for any $n\in\N^*,$ $sup_x \E_x (\int_0^{R_1}|f(X_s)|ds)^n  \le n! C(f)^n.$ In particular, $\E_{\nu}\xi^n\leq n!C(f)^n$ for any initial distribution $\nu.$
 
\end{prop}

\begin{proof}
We will first consider the case $n=2$, the general case can be obtained in the same way. Writing $\theta_s, \ s\geq 0 $ for the usual shift operator, defined on the canonical space by $X_u(\theta_s (\omega)) := X_{s + u}(\omega),$ (see \cite {RY} Chapter I, ¤ 3, p. 34) we obtain
\begin{eqnarray*}%\label{moment2}
 \left(\int_0^{R_1}|f(X_s)|ds\right)^2&=& \int_{0}^{R_{1}} \int_{0}^{R_{1}} |{f}(X_{s})| | {f}(X_{u})|ds du  \nonumber \\
& =&2! \int_{0}^{R_1}ds|f(X_s)| \int_s^{R_1}|f(X_u)|du \nonumber \\
&=&2! \int_{0}^{\infty}ds\left (|f(X_s)|\1_{\{0<s<R_1\}} \int_s^{R_1}|f(X_u)|du.\right )\nonumber \\
&\leq& 2! \int_{0}^{\infty}ds\left (|f(X_s)|\1_{\{0<s<R_1\}} \int_s^{R_1\circ\theta_s}|f(X_u)|du\right ) . 
\end{eqnarray*}
Taking expectation and using Markov's property in the last integral  gives an 
upper bound
\begin{eqnarray*}%\label{momentp2}
 \E_x\left(\int_0^{R_1}|f(X_s)|ds\right)^2 
&\leq& 2! \int_{0}^{\infty}ds\E_x\left [|f(X_s)|\1_{\{0<s<R_1\}}\E_x\left( \int_s^{R_1\circ\theta_s}|f(X_u)|du|{\F}_s\right)\right ] \\
&= &2! \int_{0}^{\infty}ds\E_x\left [|f(X_s)|\1_{\{0<s<R_1\}}\E_{X_s}\left( \int_s^{R_1}|f(X_u)|du\right)\right ]\\
&\leq& 2!C(f)^2 .
\end{eqnarray*}
Applying this argument $n$ times successively yields the  result for arbitrary $n\in\N^*.$ \end{proof}

The following estimates will also be useful in the sequel. They are obtained using local time, hence the result is
typically one-dimensional in spirit. Let $\{L_t^a ,\; t\geq 0,\; a\in\R\}$ be a local time associated to the semi-martingale $\{X_t , \; t\geq 0\},$ i.e.~a continuous increasing process such that for all $ a\in\R,$
$$|X_t-a|=|X_0-a|+\int_0^t sgn (X_s-a)dX_s +L_t^a.$$

\begin{lem}\label{lowerboundC}
Suppose that conditions \ref{model} and \ref{posrec} hold. Then for any bounded  $f:\R\to\R$, having compact support $K,$
$$ C(f) \le k \mu (|f|),$$
where $k$ is a finite constant given by 
\begin{equation*}%\label{eq:lowerboundC}
 k := \frac{M}{2}  \sup_{y\in K}s (y)\sup_{y \in K }\sup_x \E_x  L^y_{R_1} .
\end{equation*}
\end{lem}

\begin{proof}
Using assumption \ref {model},  $\sigma^2, s $ and $m$ are  continuous and strictly positive.
Using the occupation time formula, 
\begin {equation}
\sup_x\E_x\int_0^{R_1} |f(X_s)|ds\leq \int_{-\infty}^{+\infty}|f(y)|\frac1{\sigma^2(y)}\sup_x\E_xL^y_{R_1}dy
\leq \frac{M}{2}\sup_{y\in K}s(y)\sup_{y \in K }\sup_x \E_x  L^y_{R_1}\mu(|f|). 
\end{equation} 
It suffices to show that $\sup_{y \in K }\sup_x \E_x  L^y_{R_1}$ is finite. 
We start by showing that for all $ y\in\R$  $\sup_x \E_x  L_{R_1}^y=E_y  L_{R_1}^y$, which can be seen as follows:
\begin{multline*}%\label{}
 \E_x  L_{R_1}^y=\E_x  L_{R_1}^y\1_{\{R_1>T_y\}}\le  \E_x[\1_{\{R_1>T_y\}} \E_x  (L^y_{T_y+R_1\circ\theta_{T_y}}|F_{T_y})]\\
\le \P(R_1>T_y) \E_y  L_{R_1}^y\le \E_y  L_{R_1}^y.
\end{multline*}
Hence $\sup_x \E_x  L_{R_1}^y = \E_y  L_{R_1}^y$. Let $c=\inf K, d=\sup K. $ Now for $y\in K $ we write
$$\E_y  L_{R_1}^y\le \E_y  L_{R_1}^c+\E_y|L_{R_1}^y-L_{R_1}^c|\le \E_c  L_{R_1}^c+\E_y|L_{R_1}^y-L_{R_1}^c|.$$
But
\begin{equation*}%\label{eq:difloctime}
|L_{R_1}^y-L_{R_1}^c|\le |y|+|\int_0^{R_1}\1_{\{c<X_s<y\}}\sigma(X_s)dW_s|+\int_0^{R_1}\1_{\{c<X_s<d\}}|\beta (X_s)|ds .
\end{equation*}
$y\to\E_yR_1$ is continuous (see Theorem \ref{indepofpoint} below). Taking expectation with respect to $\E_y$ and taking $\sup_{y\in K }$, using continuity of $\beta$ and of $y\to\E_yR_1,$ we only need to show that 
\begin{equation*}
\sup_{y\in K }\E_y|\int_0^{R_1}\1_{\{c<X_s<y\}}\sigma(X_s)dW_s|<\infty .
\end{equation*}
By norm inclusion and isometry,
\begin{eqnarray*}
\E_y|\int_0^{R_1}1_{\{c<X_s<y\}}\sigma(X_s)dW_s|&<& \left(\E_y\left(\int_0^{R_1}\1_{\{c<X_s<y\}}\sigma(X_s)dW_s\right)^2\right)^{1/2} \nonumber \\
&\le& \left(\E_y\left(\int_0^{R_1}\1_{\{c<X_s<d\}}\sigma^2(X_s)ds\right)\right)^{1/2} .
\end{eqnarray*}
Using the continuity of $\sigma^2$ and of the map $y\mapsto\E_yR_1$ we see that $\sup_{y\in K }\E_y(\int_0^{R_1}\1_{(c<X_s<d)}\sigma^2(X_s)ds)<\infty$.
\end{proof}

We now define the point process associated to the life cycle decomposition $R_n .$ Let $N_0 = 0 $ and put for $t > 0, $
$$  N_t = \sup \{ n : R_n \le t \}=\sum _{n=1}^{\infty}\1_{\{R_n\le t\}} .$$

Then the key fact for our proof of the deviation inequality is that the processes $(N_t)_{t\geq 0}$ and $(R_n)_{n\in \N}$ are mutually inverse in the following sense :
\begin{equation*}
\{N_t\geq n\}=\{R_n\leq t\} \quad\mbox{and}\quad \{N_t\leq n\}=\{R_n\geq t\}.
\end{equation*}

\begin{lem}\label{inverses}
Suppose that $X$ verifies Assumptions \ref{model} and \ref{posrec}. Then the quantities $\E_aR_1 $ and $\E_{\mu}N_1$ are positive and finite, and
for any initial distribution $\nu$  the followings statements hold. 
\begin{enumerate}
	\item $\lim_{n\to\infty}{R_n}/n= \E_aR_1 \quad \P_{\nu}- a.s.$
	\item $\lim_{t\to\infty}{N_t}/t= \E_{\mu}N_1\quad  \P_{\nu}-a.s.$
	\item $\E_aR_1=1/{\E_{\mu}N_1}.$
\end{enumerate}
\end{lem}

\begin{proof} The finiteness of $\E_aR_1 $ follows from positive recurrence.
Statement $1.$ is the strong law of large numbers since we can write
$$\frac{R_n}{n}=\frac{R_1}{n}  +\frac1n\sum_{k=1}^{n-1}(R_{k+1}-R_{k}).$$
Using the recurrence property, $R_1<\infty \ a.s.$ and hence $R_1/n\to 0 $ almost surely. Using proposition \ref{iid} the variables $R_{k+1}-R_{k}, k\geq 1,$ are i.i.d. and equal in  law to $R_1$ under $\P_a$.
To prove the third statement we write:
$$\lim_{t\to\infty}\frac{N_t}t=\lim_{n\to\infty}\frac{N_{R_n}}{R_n}=\lim_{n\to\infty}\frac{n}{R_n}.$$
Statement $2. $ follows from the Ergodic Theorem : $(N_t)_t$ is an integrable additive functional of $X$, hence
$\lim_{t\to\infty}{N_t}/t= \E_{\mu}N_1/E_{\mu}1=\E_{\mu}N_1$ almost surely.
\end{proof}

The following proposition will be useful in the sequel:
\begin{prop}\label{muf} Suppose that $X$ verifies Assumptions \ref{model} and \ref{posrec}.
Denote $l:=\E_{\mu}(N_1).$ Then for any initial measure $\nu,$ 
$$\E_{\nu}\int_{R_1}^{R_2}f(X_s)ds= \frac{\mu(f)}{l}=\mu(f)\E_aR_1.$$
In particular, we have 
$$ |\mu (f)| \le l \cdot C(f) .$$
\end{prop}
\begin{proof}
Using the Ergodic Theorem, almost surely, 
$$\mu(f)=\lim_{t\to\infty}\frac{\int_0^tf(X_s)ds}t=\lim_{n\to\infty}\frac{\int_0^{R_n}f(X_s)ds}{R_n}.$$
On the other hand, using the strong law of large numbers,
$$\lim_{n\to\infty}\frac{\int_0^{R_n}f(X_s)ds}{R_n}=\lim_{n\to\infty}\frac{\frac1n\int_0^{R_n}f(X_s)ds}{R_n/n}=\frac{\E_{\nu}\int_{R_1}^{R_2}f(X_s)ds}{\E_aR_1}.$$
\end{proof}

\section{The deviation inequalities} 
In this section we prove the deviation inequality \eqref{erginequality}. As explained in the introduction, we use the regeneration method, which consists to ``cut the trajectory of the process into i.i.d. blocs''. However, the number of blocs before a fixed $t>0$ is a random quantity $N_t$. So in Theorem \ref{thm: Npolbound} we study the deviations of this random quantity around its mean. After that, we prove the deviation inequality for a bounded function $f$ in Theorem \ref{mainbounded}, and for $f$-bounded and compactly supported in Theorem \ref{mainint}. In the first case the dependence on $f$ is expressed through its sup-norm and in the second case through its $L^1(\mu)-$ norm.   

Throughout this section we impose assumptions \ref{model} and \ref{posrec}. Hence the measure $\mu$ of (\ref{eq:mu}) is the unique invariant probability measure of the process. 
 
\subsection{Deviations for $(N_t/t)_{t \geq 0}$}
This section is devoted to the study of deviations of $(N_t/t)_{t \geq 0}$ around its limit value $\E_{\mu}(N_1).$  The control of deviations of $(N_t/t)_{t \geq 0}$ will allow us to control the deviations of other additive functionals. We recall that $l = \E_\mu (N_1).$ The main idea of the proof of this theorem is that the processes $(N_t)$ and $(R_n)$ are mutually inverse in the sense of the lemma \ref{inverses}. The deviations of $N_t$ can therefore be expressed in terms of the deviations of $R_n=\sum_{k=0}^{n-1}(R_{k+1}-R_k)$, which is a sum of i.i.d. variables. 

\begin{theo}\label{thm: Npolbound}
Grant assumptions \ref{model} and \ref{posrec}.  Let $\nu$ be any initial distribution and $0< \ge <  1. $ Suppose that there exists $ p > 1 $  such that
 $\E_{\nu}(R_1)^{p/2}<\infty $ and $\E_{\nu} (R_2-R_1)^p<\infty.$ 
Then there exists a positive constant $C(l,p ,\nu)$ such that %for $t\geq 4/l\ge$ 
the following inequality holds:

If $p \geq 2, $ then
\begin {equation*}%\label{Npolbound}
\P_{\nu}\left(|N_t/t-l|>l\ge\right )\leq 
C(l,p, \nu)
\frac1{\ge^p}\frac1{t^{p/2}} .
\end {equation*}
If $1<p < 2 $ and $ t \geq 1,$
\begin {equation*}%\label{Npolbound}
\P_{\nu}\left(|N_t/t-l|>l\ge\right )\leq 
C(l,p, \nu)
\frac1{\ge^p}\frac1{t^{\frac{p-1}{2}}} .
\end{equation*}
Here  $C(l,p,\nu)$ is given by 
$$C(l,p,\nu) =
\left\{ 
\begin{array} {ll}
 2^{p/2}\E_{\nu}|R_1-1/l|^{p/2}+2^{3p/2}C^p_p  \E_{\nu}|\bar{\eta}_1|^p  \, l^{\frac{p}{2}}& \mbox{ if } p \geq 2\\
2^{p/2}\E_{\nu}|R_1-1/l|^{p/2}+2^{(3p+1)/2}C^p_p  \E_{\nu}|\bar{\eta}_1|^p  \, l^{\frac{p+1}{2}}&   \mbox{ if } p \in ] 1, 2 [ 
\end{array}
\right\} ,
$$
where $\bar \eta_1 = (-1)(R_2 - R_1 - \frac1l )$ and where $C_p$ is the constant of the Burkholder-Davis-Gundy inequality. 
\end{theo}

\begin{proof}
Firstly we decompose:
\begin{equation}\label{Ntfirst}
\P_{\nu}\left(|N_t/t-l|>l\ge\right )
\leq \P_{\nu}\left(N_t/t>l(1+\ge)\right )+\P_{\nu}\left(N_t/t<l(1-\ge)\right ) .
\end{equation}

Put for $k\geq 1,$
$\bar{\eta}_k=-1(R_{k+1}-R_{k}-1/l).$
For the first term of \eqref{Ntfirst}, we have 
\begin{eqnarray}\label{Ntfirst1}
\P_{\nu}\left(N_t/t>l(1+\ge)\right )&=&\P_{\nu}\left(N_t\geq [tl(1+\ge)]+1\right)=\P_{\nu}\left( R_{[tl(1+\ge)]+1}\leq t\right)\nonumber \\
&=&\P_{\nu}\left (\sum_{k=0}^{[tl(1+\ge)]}(R_{k+1}-R_{k})\leq t\right)\nonumber \\
&=&\P_{\nu}\left (\sum_{k=0}^{[tl(1+\ge)]}(R_{k+1}-R_{k} - \frac1l )\leq t\left(1-\frac{[tl(1+\ge)]+1}{tl}\right)\right) \nonumber \\
&\leq &\P_{\nu}\left (\sum_{k=0}^{[tl(1+\ge)]}(R_{k+1}-R_{k}- \frac1l)\leq t\left(1-\frac{tl(1+\ge)}{tl}\right)\right)\nonumber \\
&\leq
&\P_{\nu}\left (\sum_{k=0}^{[tl(1+\ge)]}(R_{k+1}-R_{k} - \frac1l)\leq -t\ge\right)  \nonumber\\
&\leq& \P_{\nu}\left(R_1- 1/l\leq -t\ge/2\right)+
\P_{\nu}\left(\sum_{k=1}^{[tl(1+\ge)]}\bar{\eta}_k \geq t\ge/2\right) .
\end {eqnarray}

In an analogous way,  
we treat the second term in \eqref{Ntfirst}:
\begin{eqnarray}\label{Ntfirst2}
\P_{\nu}\left(N_t/t<l(1-\ge)\right )&=&\P_{\nu}\left (N_t\leq [tl(1-\ge)]\right )= \P_{\nu}\left(R_{[tl(1-\ge)]}\geq t\right ) \nonumber \\
&=&\P_{\nu}\left (\sum_{k=0}^{[tl(1-\ge)]-1}(R_{k+1}-R_{k}-\frac{1}{l})\geq t \left( 1 -\frac{[tl(1-\ge)]}{t l} \right) \right)\nonumber \\
 &\leq& \P_{\nu}\left(R_1-\frac1l\geq t\ge/2\right)
+\P_{\nu}\left(\sum_{k=1}^{[tl(1-\ge)]-1}\bar{\eta}_k\leq -t\ge/2\right).
\end{eqnarray}
Let $M_0=0$ and $M_n=\sum_{k=1}^{n}\bar{\eta}_k.$ For $k\geq 1,$ $(\bar{\eta}_k)$ are i.i.d. centered random variables such that $\E_{\nu}|\bar{\eta}_k|^p <\infty,$ hence $(M_n)_{n \geq 1}$ is an $L^p$ martingale such that $[M]_n=\sum_{k=1}^n\bar{\eta}_k^2.$ Denote $M^*_n=\sup_{k\leq {n}}|M_k|.$ 
As a consequence of \eqref{Ntfirst1} and \eqref{Ntfirst2} we can write
\begin{equation}\label{Ntsecond}
\P_{\nu}\left(|N_t/t-l|>l\ge\right )
\leq \P_{\nu}\left(|R_1-1/l|\geq t\ge/2\right )+\P_{\nu}\left( M^*_{[tl(1+\ge)]}\geq t\ge/2\right) .
\end{equation}

We use the Burkholder-Davis-Gundy inequality to bound the last term in \eqref {Ntsecond}.   
By the Burkholder-Davis-Gundy inequality, for all $p>1$ there exists a constant $C_p$ depending only $p$ such that $\|M_n^*\|_p\leq C_p\|[M]_n^{1/2}\|_p,$  hence
%\begin{equation}\label{BGD}
 $\E_{\nu}(M_n^*)^p\leq C_p^p\E_{\nu}\left (\sum_{k=1}^n\bar{\eta}_k^2\right)^{p/2}.$
%\end{equation}

If $p\geq 2$, using H\"older's inequality, 
\begin{equation}\label{BGD1}
\left(\sum_{k=1}^n\bar{\eta}_k^2\right)^{p/2}\leq n^{\frac{p}{2}-1}
\sum_{k=1}^n|\bar{\eta}_k|^p ,\quad\mbox{hence}\quad \E_{\nu} (M_n^*)^p\leq C_p^p n^{p/2}\E|\bar{\eta}_1|^p.
\end{equation}
If $1<p<2,$ using H\"older's inequality together with the sub-additivity of the function $x \mapsto \sqrt{x},$ 
\begin{equation}\label{BGD2}
  \left(\sum_{k=1}^{n}\bar{\eta}_k^2\right)^{p/2 }\le n^{\frac{p-1}{2}} \sum_{k=1}^{n}|\bar{\eta}_k|^p ,\quad\mbox{hence}\quad \E_{\nu} (M_n^*)^p\leq C_p^p n^{\frac{p+1}2}\E|\bar{\eta}_1|^p.
  \end{equation}
 Finally, if $p\geq 2,$
\begin {eqnarray*}%\label{NNpolbound}
\P_{\nu}\left(|N_t/t-l|>l\ge\right )&& 
\le   \frac{2^{p/2}\E_{\nu}|R_1-1/l|^{p/2}}{(t\ge)^{p/2}}+
2^p C^p_p  \E_{\nu}|\bar{\eta}_1|^p \, [tl(1+\ge)]^{p/2}\frac{1}{(\ge t)^{{p}}} \\
&& \leq \left( 2^{p/2}\E_{\nu}|R_1-1/l|^{p/2}+2^{3p/2}C^p_p  \E_{\nu}|\bar{\eta}_1|^p  \, l^{\frac{p}{2}} \right) \frac1{\ge^p}\frac{1}{t^{\frac{p}{2}}},
\end {eqnarray*}
and if $1<p<2,$ for $t \geq 1,$
\begin {eqnarray*}%\label{NNpolbound}
\P_{\nu}\left(|N_t/t-l|>l\ge\right ) 
&& \le   \frac{2^{p/2}\E_{\nu}|R_1-1/l|^{p/2}}{(t\ge)^{p/2}}+
2^p C^p_p  \E_{\nu}|\bar{\eta}_1|^p \,[tl(1+\ge)]^{(p+1)/2}\frac{1}{(t\ge)^{p}} 
\\
&& \leq \left( 2^{p/2}\E_{\nu}|R_1-1/l|^{p/2}+2^{(3p+1)/2}C^p_p  \E_{\nu}|\bar{\eta}_1|^p  \, l^{\frac{p+1}{2}} \right) \frac1{\ge^p}\frac{1}{t^{\frac{p-1}{2}}}.
\end {eqnarray*}

\end{proof}
\subsection{Rate of convergence in the Ergodic Theorem} 
We apply the results of the previous section to get a bound on the rate of convergence in the Ergodic Theorem for additive functionals $\int_0^t f(X_s)ds$, where $f\in\L^1(\mu) .$ We consider two situations. Firstly, the case where $f$ is bounded, secondly, the case where $f$ is bounded and compactly supported.
Our bound depends on $f$ through $\|f\|_\infty$ in the first case, and through $\mu(|f|)$ in the second one. %Actually, in the second case  we only need the finiteness of $\mu(f)$ and $C(f)$. Such conditions are often used in the study of recurrent diffusions. Functions $f$ with finite $C(f)$ are called special with respect to the process $X$ and to the decomposition into life cycles $(R_n)_n$. It is known that for strong Feller diffusions all bounded functions with compact support are special with respect to all decompositions. That's why we work in the second case with this class of functions.  
% Recall that we defined $ C(f) = \sup_x \E_x \int_0^{R_1} |f|(X_s) ds $ in the last case. 
In both proofs we use the following decomposition of trajectories: the trajectory before $R_1$, the trajectory between $R_k$ and $R_{k+1}$, $1\leq k\leq N_t+1$, and finally the trajectory between $t$ and $N_t+1.$ We also restrict this decomposition to the set $\Omega_t$ where $N_t$ is close to its mean. Hence the main term -- the sum of parts between $R_k$ and $R_{k+1}$ -- becomes just a sum of i.i.d. variables. The control of the complementary of $\Omega_t$ is given by the theorem \ref{thm: Npolbound}.  
 
\begin{theo}\label{mainbounded}
Grant assumptions \ref{model} and \ref{posrec}. Let $f\in\L^1(\mu)$. Suppose that $\|f\|_\infty <\infty$.
 Let $\nu$ be any initial distribution and $0< \ge <  \|f\|_\infty  . $ Suppose that there exists $p>1$  such that
 $\E_{\nu}(R_1)^{p/2}<\infty $ and $\E_{\nu}(R_2-R_1)^p<\infty.$  
%\begin{description} 
%\item[ (a)] 
Then for all $t\geq 1$  the following inequality holds:
\begin{equation*}%\label{polbounded}
\P_{\nu}\left(\left|\frac 1t\int_0^tf(X_s)ds-\mu(f)\right|>\ge\right)\leq 
\left\{
\begin{array}{ll}
K(l,p, \nu,X)\frac 1{\ge^p}\|f\|_\infty ^p \, {t^{- p/2}} & \mbox{ if } p \geq 2 \\
K(l,p, \nu,X)\frac 1{\ge^p}\|f\|_\infty ^p \, {t^{-\frac{p-1}{2}}} & \mbox{ if }1< p < 2
\end{array}
\right\}. 
\end{equation*}

%\item [(b)] For a bounded function $f$ with compact support we have \begin{equation}\label{polintegr}
%\P_{\nu}\left(\left|\frac 1t\int_0^tf(X_s)ds-\mu(f)\right|>\ge\right)\leq
%\left\{  
%\begin{array}{ll}
%K(l,p, \nu )\frac 1{\ge^p}(C(f)\vee |\mu(f)|)^p \, {t^{-p/2}}& \mbox{ if } p \geq 2\\
%K(l,p, \nu ) \frac 1{\ge^p}(C(f)\vee |\mu(f)|)^p \, {t^{- \frac{p-1}{2}}} & \mbox{ if } p < 2
%\end{array}
%\right\} , 
%\end{equation}
%\end{description}
Here $K(l, p, \nu,X)$ is a positive constant, different in the two cases, which depends on $l,p ,\nu$
and on the process $X$ through the life cycle decomposition, but which does not depend on $f$, $t$, $\ge.$
\end{theo}
\begin{rem}\label{conditions}
Since $\E_{\nu}(R_1)^p\leq 2^{p-1}(\E_{\nu}T_b^p+\E_bT_a^p)$, we can see that the hypotheses of the 
theorem \ref{mainbounded} are satisfied if $\E_aT_b^p<\infty,$ $\E_bT_a^p<\infty$ and $\E_{\nu}T_b^{p/2}<\infty.$ 
%and $\E_{\nu} T_b^{p/2}<\infty.$
%In particular, if for one couple of points $x,y$ $\E_xT_y^p<\infty$,
%then for all couple of points  $x,y$ $\E_xT_y^p<\infty$ (see the theorem \ref{indepofpoint}),
%and the theorem \ref{mainbounded} is satisfied for any $x\in\R$ and any initial measure of the form $\nu=\delta_{x}.$ 
Corollary \ref{expbounded} gives some  explicit conditions for that in terms of the coefficients of $X$.  
 \end{rem}
\begin{rem}
Using regeneration techniques developed in \cite{dashaeva} it should be possible to get some multidimensional version of the previous theorem, but on this stage we are not able to state any practical condition ensuring the existence of moments of regeneration times in this case. For that  reason in the present paper we restrict our attention to the one-dimensional diffusions.  
\end{rem}
\begin{proof}
Put $\bar f := f - \mu (f) .$ 
Recall that $0<\ge < \|f\|_\infty.$ 
Denote $\delta={\varepsilon }/{\|f\|_\infty}$ and
\begin{equation*}%\label{Omega}
\Omega_t=\left\{ \left|\frac{N_t}t-l\right|\leq l\delta \right\} .
\end{equation*}

We shall use the following decomposition.
 \begin{eqnarray*}%\label{fdecomp}
&& \P_{\nu}\left(\left|\int_0^tf(X_s)ds-t\mu(f)\right|>t \varepsilon  \right)
\nonumber\\
&& \leq\P_{\nu}\left(\left|\int_0^t \bar f(X_s)ds \right|>t \varepsilon  \ ; \Omega_t\right)+\P_{\nu}\left(\Omega_t^c\right) \nonumber \\
&& \leq\P_{\nu}\left(\left|\int_0^{R_1}\bar f(X_s)ds\right|>\frac{t\varepsilon}3\right)
+\P_{\nu}\left(\left|\int_{R_1}^{R_{N_{t}+1}} \bar f(X_s)ds \right|>\frac{t\varepsilon}3\ ; \Omega_t\right)\nonumber \\
 &&+\P_{\nu}\left(\left|\int_t^{R_{N_{t}+1}}\bar f(X_s)ds\right|>\frac{t \varepsilon}{3}; \Omega_t\right)
+\P_{\nu}\left(\Omega_t^c\right)\nonumber \\
&&=A+B+C+D . 
\end{eqnarray*}

%In the sequel, $K(p)= K(l,p, \nu)$ means a constant that may change from line to line but that does not depend on $t,\ge,f$. This constant depends only on the process, the choice of life cycles and $p$.
For the term $A,$ we have, since $ \| \bar f \|_\infty \le 2\|f\|_\infty,$  
\begin{equation}\label{A1pol}
\P_{\nu}\left(\left|\int_0^{R_1}\bar f(X_s)ds\right|>\frac{t \varepsilon}3 \right)\leq \P_{\nu}\left(R_1>\frac{t \varepsilon}{6 \|f\|_\infty}\right) \leq \frac{\E_{\nu}R_1^{p/2}}{ t^{p/2}}\left(\frac{ 6\|f\|_\infty}{ \varepsilon}\right)^{p/2} .
\end{equation}
Recall that for $n\geq 1 ,$  $\xi_n=\int_{R_n}^{R_{n+1}}\bar  f(X_s)ds$ are i.i.d.~random variables. Using proposition \ref{iid}, the law of $\xi_n$, $n\geq 1,$ does not depend on the initial distribution and is equal to the law of $\xi_0=\int_0^{R_1} \bar f(X_s)ds$ under $\P_a$. 
Recall (proposition \ref{muf}) that $ \E_a\xi_0=\mu(\bar f)/l = 0$.

In the sequel we need $\E_{\nu}|\xi_k |^p<\infty$, which can be seen as follows :  $$\E_{\nu}|\xi_k |^p \le 2^p\|f\|_\infty^p\E_{\nu}(R_2-R_1)^p .$$
 
Now we treat the term $B$, which is the main term of the decomposition. Denote $M_0 = 0,$ $M_n = \sum_{k=1}^n ( \xi_k  )$ and $M^*_n=\sup_{k=0,\ldots n}|M_k|.$ Then we have
\begin{eqnarray*}%\label{Bpol}
B&=&\P_{\nu}\left(|\int_{R_1}^{R_{N_{t+1}}} \bar f(X_s)ds |>\frac{ t \varepsilon}3;\ \Omega_t\right)\nonumber \\
%&=&\P_{\nu}\left(|\sum_{k=1}^{N_t}\xi_k-N_t\mu(f)/l|>\frac{t\ge\mu(f)}4; \Omega_t \right)\nonumber \\
&\leq& \P_{\nu}\left(|\sum_{k=1}^{N_t}\xi_k |>\frac{t \varepsilon }3 ;\quad |{N_t}/t-l|\leq l \delta \right)\nonumber \\
&\leq &\P_{\nu}\left(\sup_{n\leq [tl(1+ \delta)]}|M_n|>\frac{ \varepsilon}3 \right) 
\le \frac{3^p \E_\nu  ( \, M^*_{[tl(1+\delta)] })^{p}}{ \, \varepsilon^p \,  t^p}.
%4^p[tl(1+\ge /4)]^{p/2}\E_{\nu}|\xi_1-\E_x \xi_1|^p}{\mu(f)^p\ge^pt^p} ,
\end{eqnarray*}  
We want to use the Burkholder-Davis-Gundy inequality for the martingale $M_n.$ Now as in \eqref{BGD1}, \eqref{BGD2} we have
for $p \geq 2$ 
\begin{equation*}%\label{BGD4}
\E_{\nu}(M^*_n)^p\le  C_p^p n^{\frac{p}{2}- 1}\E_{\nu} \sum_{k=1}^{n}|\xi_k  |^p=C_p^pn^{p/2 }\E_{\nu}|\xi_1  |^p ,
\end{equation*}
and for $1< p < 2,$
\begin{equation*}%\label{BGD3}
\E_{\nu}(M^*_n)^p\le  C_p^p n^{\frac{p-1}{2}}\E_{\nu} \sum_{k=2}^{n+1}|\xi_k  |^p=C_p^pn^{\frac{p+1}{2}}\E_{\nu}|\xi_1  |^p.
\end{equation*}
Finally 
 we have for $p \geq 2,$
\begin{equation}\label{B1pol}
B
\leq \frac{C_p^p 3^p[tl(1+\delta)]^{p/2}}{ \varepsilon^p t^p}2^p\|f\|_{\infty}^p \E_{\nu}|R_2-R_1|^p\leq
K(p)\|f\|_{\infty}^p\frac1{t^{p/2}}\frac{1}{\varepsilon^p},
\end{equation}  
where $K(p)=C_p^p12^pl^{p/2}\E_{\nu}|R_2-R_1|^p ,$ and for $p < 2,$ 
\begin{equation*}%\label{B1pol2}
B \leq \frac{C_p^p3^p[tl(1+\delta)]^{\frac{p+1}{2}}}{\varepsilon^p t^p}2^p\|f\|_{\infty}^p \E_{\nu}|R_2-R_1|^p\leq (2l)^{1/2} \, K(p)\|f\|_{\infty}^p\frac1{t^{\frac{p-1}{2}}}\frac{1}{\varepsilon^p}.
\end{equation*}

For the term $C$  we can write
\begin{eqnarray*}%\label{Cpol}
C&&=\P_{\nu}\left(\left|\int_t^{R_{N_{t}+1}}\bar f (X_s)ds\right|>\frac{t \varepsilon}{3}; \Omega_t\right)
\\
&&\leq \sum_{k=1}^{[tl(1+\delta)]}\P_{\nu}\left(|\int_t^{R_{N_{t}+1}}\bar f(X_s)ds|>\frac{t \varepsilon}{3};\ N_t=k \right)
\\ 
&&\leq \sum_{k=1}^{[tl(1+\delta )]}\P_{\nu}\left(\int_{R_k}^{R_{k+1}}|\bar f|(X_s)ds>\frac{t \varepsilon }{3} \right)\\
&&\leq  tl(1+\delta )\frac{\E_{\nu}(\int_{R_1}^{R_{2}}|\bar f|(X_s)ds)^p}{t^p}\left(\frac{3}{\varepsilon }\right)^p \\
&&\leq \frac{1}{t^{p-1} } \frac{2^p \|f\|_{\infty}^p}{\varepsilon^p}2l3^p\E|R_2-R_1|^p .
%\leq \frac 1{t^{p/2}}\frac{l2\E_aR_1^p(4M)^p}{(\ge\mu(f)^p}.
\end{eqnarray*} Finally,
\begin{equation}\label{C1pol}
C\leq K(p)\frac 1{t^{p- 1}}\frac{\|f\|_{\infty}^p}{\varepsilon^p} 
\le \left\{ 
\begin{array}{ll}
K(p) \frac{1}{t^{p/2} } \frac{\|f\|_{\infty}^p}{\varepsilon^p}  & \mbox{ if } p \geq 2\\
K(p) \frac{1}{t^{p-1} } \frac{\|f\|_{\infty}^p}{\varepsilon^p}  & \mbox{ if } p <  2
\end{array}
\right\},
\end{equation}
where $K(p)=2^{p+1} l 3^p\E_aR_1^p .$

For the term $D,$ we use the theorem \ref{thm: Npolbound}  : 
$$ D \le  \left\{ 
\begin{array}{ll}
 C(p) {t^{-p/2}}\frac{\|f\|_{\infty}^p}{\varepsilon^p} & \mbox{ if } p \geq 2\\
 C(p) {t^{- \frac{p-1}{2}}}\frac{\|f\|_{\infty}^p}{\varepsilon^p} & \mbox{ if } p < 2
 \end{array}
\right\} .$$ 
Here, $C(p)$ is the constant of the theorem \ref{thm: Npolbound}.
Finally we obtain, putting together \eqref{A1pol}, \eqref{B1pol} and \eqref{C1pol} 
\begin{equation*}%\label{F1pol}
\P_{\nu}\left (|\int_0^tf(X_s)ds-t\mu(f)|>t  \varepsilon \right)\leq K(l,p,\nu,X)\frac{1}{t^{\alpha}}\left (\frac{\|f\|_{\infty}}{\varepsilon}\right)^p\vee \left (\frac{\|f\|_{\infty}}{ \varepsilon}\right)^{p/2},
\end{equation*}
where $\alpha=p/2$ if $p\geq 2,$ and  $\alpha=(p-1)/2$  for $1<p<2.$

Finally, since $\|f\|_{\infty} > \varepsilon , $ 
$$ \left (\frac{\|f\|_{\infty}}{\varepsilon}\right)^p\vee \left (\frac{\|f\|_{\infty}}{\varepsilon }\right)^{p/2} = \left (\frac{\|f\|_{\infty}}{\varepsilon }\right)^p .$$ 
Then the theorem  follows.
\end{proof}

In the case where $f$ is bounded and compactly supported we get the  version of the deviation inequality  with $\L^1(\mu)$ norm of $f$ instead of its sup-norm.  In some practical situations this can be of major importance. In the next theorem we only deal with integer $p \in \N^* .$ This is due to the fact that proposition \ref{prop:momentp} is only stated for integer moments.
\begin{theo}\label{mainint}
Grant assumptions \ref{model} and \ref{posrec}. Let $f $ be a bounded function with compact support. Let $\nu$ be any initial distribution. Suppose that there exists $ p\in\N ,\  p > 1,$  such that
 $\E_{\nu}(R_1)^{p}<\infty $  and $\E_{\nu}(R_2-R_1)^p<\infty.$ 
 Then for all $t \geq 1,$ for all $0<\ge<\mu(|f|)$ the following inequality holds:
\begin{equation*}%\label{polintegr}
\P_{\nu}\left(\left|\frac 1t\int_0^tf(X_s)ds-\mu(f)\right|>\ge\right)\leq K(l,p,  X )\frac 1{\ge^p}\mu(|f|)^p \, {t^{-p/2}}.
\end{equation*}
Here $K(l, p , X)$ is a positive constant which depends on $l,p, X ,$ but which does not depend on $f$, $t$, $\ge.$
\end{theo}
\begin{rem} The corollary \ref{explintegrable} gives
some explicit conditions for the theorem \ref{mainint} in terms of coefficients of $X$.  
\end{rem} 
%\begin{rem}
%Recall that $C(f)\leq k \cdot \mu(|f|) $ (see lemma \ref {lowerboundC}), hence the bound in Theorem \ref{mainint} can also be expressed in terms of $\mu(|f|).$  
%\end{rem}
\begin{proof}
Since $0<\ge < \mu (|f|),  $ we can write 
$\varepsilon =  \mu (|f|)   \delta, $ where $0 < \delta  < 1 .$ 
%We start by deriving a deviation bound for
%\begin{equation}
%\P_{\nu}\left(\left|\int_0^tf(X_s)ds-t\mu(f)\right|>t \mu(f) \delta \right) .
%\end{equation}

Denote 
\begin{equation*}%\label{Omega}
\Omega_t=\left\{ \left|\frac{N_t}t-l\right|\leq \frac{l\delta}{4} \right\} .
\end{equation*}

We shall use the following decomposition.
 \begin{eqnarray*}%\label{fdecomp2}
&& \P_{\nu}\left(\left|\int_0^tf(X_s)ds-t\mu(f)\right|>t \varepsilon \right)
\nonumber\\
&& \leq\P_{\nu}\left(\left|\int_0^tf(X_s)ds-t\mu(f)\right|>t  \ge \ ; \Omega_t\right)+\P_{\nu}\left(\Omega_t^c\right) \nonumber \\
&& \leq\P_{\nu}\left(\left|\int_0^{R_1}f(X_s)ds\right|>\frac{t\ge}4\right)
+\P_{\nu}\left(\left|\int_{R_1}^{R_{N_{t}+1}} f(X_s)ds-N_t\frac{\mu(f)}{l}\right|>\frac{t\ge}4\ ; \Omega_t\right)\nonumber \\
 &&+\P_{\nu}\left( |N_t\frac{\mu(f)}{l}-t\mu(f)|>\frac{t\ge}4;\ \Omega_t\right)+\P_{\nu}\left(\left|\int_t^{R_{N_{t}+1}}f(X_s)ds\right|>\frac{t\ge}{4}; \Omega_t\right)
+\P_{\nu}\left(\Omega_t^c\right)\nonumber \\
&&=A+B+E+C+D . 
\end{eqnarray*}
All the long of the proof $K$ is a positive constant, not always the same, which depends on $l,p $
and on the process $X$ through the life cycle decomposition, but which does not depend on $f$, $t$, $\ge.$

We start with the term $E.$ Using ${\mu(|f|)}/{|\mu(f)|}\geq 1$ together with $ x/0=+\infty$ for $x>0,$ we have
\begin{equation}\label{E}
E=\P_{\nu}\left( \left|\frac{N_t}{t}-l\right|>\frac{\mu(|f|)}{|\mu(f)|}\frac{l\delta}4;\ \Omega_t\right)=0 .
\end{equation}

For the term $A,$ we have, applying proposition \ref{prop:momentp} and lemma \ref{lowerboundC},
\begin{multline}\label{A2pol}
A=\P_{\nu}\left(\left|\int_0^{R_1}f(X_s)ds\right|>\frac{t \ge}4\right) \leq\frac{ \E_{\nu}(\int_0^{R_1}|f(X_s)|ds)^{p}}{ t^{p}}\left(\frac{4}{\ge }\right)^{p} 
\leq\frac{ {p!}\, C(f)^{p}}{t^{p/2}}\left(\frac{4}{\ge}\right)^{p}\leq\\
\leq k^p  \frac{ {p!}\, \mu(|f|)^{p}}{t^{p/2}}\left(\frac{4}{\ge}\right)^{p} .
\end{multline} Recall that for $n\geq 1 ,$  $\xi_n=\int_{R_n}^{R_{n+1}}f(X_s)ds$ are i.i.d. equal in law to the
$\xi_0$ under $\P_a$. By proposition \ref{muf} $ \E_a\xi_0=\mu(f)/l$. Write $M_0 = 0,$ $M_n = \sum_{k=1}^n ( \xi_k - E_{\nu} (\xi_k) ).$ 
 We have that 
$\E_{\nu}|\xi_k-\mu(f)/l|^p<\infty$, since  
$$\E_{\nu}|\xi_k-\mu(f)/l|^p<2^p\left(\E_{\nu}|\xi_k|^p+|\mu(f)/l|^p\right)< 2^pp!C(f)^p+2^p|\mu(f)/l|^p<\infty. $$

Then
\begin{eqnarray*}%\label{Bpol}
B=\P\left(|M_{N_t}|>\frac{t\ge}4;\  \Omega_t \right )&\leq &\P_{\nu}\left(\sup_{n\leq [tl(1+\delta /4)]}|M_n|>\frac{t\ge}4 \right) 
\le \frac{ 4^p\E_\nu  ( \, M^*_{[tl(1+\delta /4)] })^{p}}{\ge^p \,  t^p}.
%4^p[tl(1+\ge /4)]^{p/2}\E_{\nu}|\xi_1-\E_x \xi_1|^p}{\mu(f)^p\ge^pt^p} ,
\end{eqnarray*}  
As in the proof of Theorem \ref {mainbounded}, we use the Burkholder-Davis-Gundy inequality for the martingale $M_n.$ Now as in the proof of 
\eqref{BGD1}, since $p \geq 2,$  
\begin{equation*}%\label{BGDbis3}
\E_{\nu}(M^*_n)^p\le C_p^p n^{p/2 }\E_{\nu}|\xi_1 - \E_x ( \xi_1) |^p\leq C_p^pn^{p/2 }\left (2^pp!C(f)^p+2^p|\mu(f)/l|^p \right ). 
\end{equation*}
Hence, since $ C (f) \le k \mu (|f|)  $  (Lemma \ref{lowerboundC}),
\begin{equation}\label{B2pol}
B
\leq \frac{K(p)[tl(1+\delta /4)]^{p/2}}{\ge^pt^p} \mu (|f|)^p
\leq
K(l,p, X)\mu(|f|)^p\frac1{t^{p/2}}\frac{1}{\ge^p}.
\end{equation} 
%where $K(p)=C_p^p16^pl^{p/2}(p!+1 ).$

For the term $C$ as in the proof of Theorem \ref{mainbounded} and using Proposition \ref{prop:momentp}  we can write:
\begin{multline*}%\label{Cpol}
C=\P_{\nu}\left(\left|\int_t^{R_{N_{t}+1}}f(X_s)ds\right|>\frac{t \ge}{4}; \Omega_t\right)
\leq  tl(1+\delta /4)\frac{\E_{\nu}(\int_{R_k}^{R_{k+1}}|f|(X_s)ds)^p}{t^p}\left(\frac{4}{\ge}\right)^p\\
\leq   tl(1+\delta /4)\frac{p!C(f)^p}{t^p}\left(\frac{4}{\ge}\right)^p.
%\leq \frac 1{t^{p/2}}\frac{l2\E_aR_1^p(4M)^p}{(\ge\mu(f)^p}.
\end{multline*}
We get 
\begin{equation}\label{C2pol}
C\leq K(p)\frac 1{t^{p- 1}}\frac{C(f)^p}{\ge^p} 
\le 
K( p) \frac{1}{t^{p/2} } \frac{C(f)^p}{\ge^p} 
\le 
K(l,p,X) \frac{1}{t^{p/2} } \frac{\mu(|f|)^p}{\ge^p} 
\end{equation}
since $ p \geq 2 ,$ using once more Lemma \ref{lowerboundC}. 
%where $K(p)=2l4^pp!$

For the term $D$ we have : 
\begin{equation}\label{D}
 D=\P_{\nu}\left ( |\frac{N_t}t-l|\geq\frac{l\delta}4\right )
  \le  K  {t^{-p/2}}\frac1{\delta ^p} \le  K(l,p )  {t^{-p/2}}\frac{\mu(|f|)^p}{(\delta \mu(|f|))^p} = 
K(l,p) \frac{1}{t^{p/2} } \frac{\mu(|f|)^p}{\ge^p}  .
\end{equation}
Finally, we put together
 \eqref{E}, \eqref{A2pol}, \eqref{B2pol}, \eqref{C2pol} and \eqref{D}, and the theorem  follows.
\end{proof}

\section{ Kac formula }

In Theorems \ref{thm: Npolbound}, \ref{mainbounded} and \ref{mainint}, the speed of convergence is governed by the $p$-th moment of the regeneration time, which can be expressed in terms of  $\E_xT_y^p.$  In this section we give a generalized version of Kac's moments formula (compare to \cite{FitsKac} and \cite{Val}). It will be used to prove that the moments $\E_xT_y^p,$ $p\geq 1,$ exist (or not) simultaneously for all couples $x<y$, (resp. $x>y$), see the theorem \ref{indepofpoint}. Also, Kac's formula will be used in the last section to give necessary and sufficient conditions of existence of such a moments. 

Fix any pair of points $a,b$ with  $- \infty < a < b < + \infty .$ For  $a\leq x\leq b$ let us consider
\begin{equation*}%\label{Tab}
T_{a,b}=\inf \{t\geq 0; \ X_t\notin ]a,b[\} . 
\end{equation*}

Let $G$ be the Green's function %potential kernel
associated to the stopping time $T_{a,b},$ defined by
\begin{equation*}%\label{GGab}
G(a,b,x,\xi)=\left\{\begin{array}{lr}
\frac{(S(x)-S(a))(S(b)-S(\xi))}{S(b)-S(a)}& a\leq x\leq \xi\leq b\\
\frac{(S(b)-S(x))(S(\xi)-S(a))}{S(b)-S(a)}& a\leq \xi \leq x\leq b\\
0 & \mbox{otherwise.} 
\end{array}\right.
\end{equation*}
\begin{theo}\label{thm:mainKac}(Generalized Kac's moment formula.)

Let $f:\R\to\R$ be such that the function $x\to \E_xf'(T_{a,b})$ is continuous on $[a,b].$ Then
\begin{equation}\label{mainKac}
\E_xf(T_{a,b})=\int_{-\infty}^{+\infty} G(a,b,x,\xi)\E_{\xi}f'(T_{a,b})m(\xi)d\xi . 
\end{equation}
\end{theo}
 \begin{proof}
For any $f:\R^+\to\R$ denote $u_{f}(x)=u_{f}(x,a,b)=\E_{x}f(T_{a,b})$ and
 let $$u(x)=u(x,a,b)=\int_{-\infty}^{+\infty} G(a,b,x,\xi)u_{f'}({\xi})m(\xi)d\xi.$$
 We see that $u$ is continuous on $[a,b].$ Let $(L u) (x)=\frac12 \sigma^2(x) u'' (x) + \beta (x) u' (x)  $ be the generator  of the semi-group of $X.$
 An easy calculation using the derivation of an integral with variable upper limit and $LS=0$ shows that under our assumption $ u$ satisfies
 \begin{equation*}%\label{LL}
\left\{\begin{array}{ll}
Lu(x)=-u_{f'}(x),& a< x< b  \\
u(a)=a(b)=0 .
\end{array}\right. 
\end{equation*}

Hence
the Ito formula applied to  $u $  gives
\begin{equation*}%\label{ddu1}
du(X_t)=-u_{f'}(X_t)dt+dM_t;\; u(X_t)=u(x)-\int_0^tu_{f'}(X_s)ds+M_t , 
\end{equation*}
where $M_t=\int_0^tu'(X_s)\sigma^2(X_s)dW_s$ is a continuous local martingale such that $M_{t\wedge T_{a,b}}$ is uniformly integrable.
Doob's stopping rule gives
$$0=u(X_{T_{a,b}})=u(x)-\int_0^{T_{a,b}}u_{f'}(X_s)ds+M_{T_{a,b}},$$
thus
$$u(x)=\E_x\int_0^{T_{a,b}}u_{f'}(X_s)ds . $$
Then
\begin{multline*}
u(x)=\E_x\int_0^{T_{a,b}}\E_{X_s}f'(T_{a,b})ds=\int_0^{\infty}
\E_x(1_{\{s<T_{a,b}\}}\E_x(f'(T_{a,b}\circ\theta_s)|{\cal F}_s)ds)=\\
=\int_0^{\infty}\E_x(\E_x(f'(T_{a,b}\circ\theta_s)1_{\{s<T_{a,b}\}}|{\cal F}_s)ds)
=\int_0^{\infty}\E_x(\E_x(f'(T_{a,b}-s)1_{\{s<T_{a,b}\}}|{\cal F}_s)ds)=\\
=\int_0^{\infty}
\E_x(f'(T_{a,b}-s)1_{\{s<T_{a,b}\}})ds=\E_x\int_0^{T_{a,b}}f'(T_{a,b}-s)ds= \E_xf(T_{a,b}),
\end{multline*}
and the theorem follows.
 \end{proof}

Define for $x<b$
\begin{equation*}%\label{Gb}
G(-\infty,b,x,\xi)=\left\{\begin{array}{rr}
(S(b)-S(x)),&-\infty< \xi \leq x\\
(S(b)-S(\xi)),&  x\leq \xi\leq b\\
0,& \xi\geq b
\end{array}\right.
\end{equation*}
and for $x>a,$
\begin{equation*}%\label{Ga}
G(a,+\infty,x,\xi)=\left\{\begin{array}{rr}
0,& \xi\leq a\\
(S(\xi)-S(a)),&  a\leq \xi\leq x\\
(S(x)-S(a)),& x \leq \xi< \infty .
 \end{array}\right.
\end{equation*}

\begin{prop}\label{prop:kac}
Under assumptions \ref{model} and \ref{posrec}, we have for all $p \in [1, + \infty [,$ $a\in\R$, $b\in\R$,  
\begin{equation}\label{Kacb} 
\E_xT_b^p=p  \int_{-\infty}^{+\infty}G(-\infty,b,x,\xi)\E_{\xi} T_{b}^{p-1} m(\xi)d\xi\quad  \forall x<b 
\end{equation}
and 
\begin{equation}\label{Kaca} 
\E_xT_a^p =p \int_{-\infty}^{+\infty}G(a,+\infty,x,\xi)\E_{\xi} T_{a}^{p-1} m(\xi)d\xi\quad \forall x>a .
\end{equation}
\end{prop}

\begin {rem}
The expressions \eqref{Kacb} and \eqref{Kaca} are always defined, because all functions we integrate are positive. In Theorem \ref{indepofpoint} below we discuss the issue of finiteness of these terms. %Note also that $s, m $ and $S$ are continuous. We even have $S\in C^2(\RR)$ and $LS=0$.
\end{rem}

\begin{proof}
As an application of \eqref{mainKac} with $f(x)=x$ we obtain
\begin{equation}\label{kac1} 
\E_xT_{a,b}=\int_{-\infty}^{+\infty}G(a,b,x,\xi)m(\xi)d\xi.
\end{equation}
As a consequence, being an integral with variable upper limit and continuous integrand, the function $x\mapsto \E_xT_{a,b}$ is continuous on $[a,b].$ Thus \eqref{mainKac} with $f(x)=x^2$ applies and gives
$$\E_xT^2_{a,b}=2\int_{-\infty}^{+\infty}G(a,b,x,\xi)\E_{\xi}T_{a,b}m(\xi)d\xi,$$
which is also  continuous being an integral with variable upper limit and continuous integrand. Finally, after $n$ applications of \eqref{mainKac} we have
\begin{equation}\label{intKac}
\E_xT^n_{a,b}=n\int_{-\infty}^{+\infty}G(a,b,x,\xi)\E_{\xi}T_{a,b}^{n-1}m(\xi)d\xi.
\end{equation}
Using monotone convergence, we get
$$\E_xT_b^n=\lim_{a\to -\infty}\E_xT_{a,b}^n.$$
Note that
\begin{equation*}%\label{GGb}
G(-\infty,b,x,\xi)=\lim_{a\to-\infty}G(a,b,x,\xi)
=\left\{\begin{array}{rr}
(S(b)-S(\xi))&  x\leq \xi\leq b\\
(S(b)-S(x))& \xi \leq x\leq b \\
0& \xi >b. 
\end{array}\right.
\end{equation*}
Moreover, for all $ a<x<b,$ $G(a,b,x,\xi)\leq G(-\infty,b,x,\xi).$
So, 
if the integral 
\begin{equation}\label{eq:integral}
 \int_{-\infty}^{+\infty}G(- \infty ,b,x,\xi)\E_{\xi}T_{a,b}^{n-1} m(\xi)d\xi
 \end{equation}
converges, using dominated convergence,  we pass to the limit when $a\to\-\infty$, which gives 

\begin{equation*}%\label{EExTbn}
\E_xT_{b}^n=n\int_{-\infty}^{+\infty}G(-\infty,b,x,\xi) \E_{\xi}T_{b}^{n-1} m(\xi)d\xi.
\end{equation*}
 
We can rewrite this last expression as
\begin{equation*}%\label{EExTbn}
\E_xT_{b}^n=n\left((S(b)-S(x))\int_{-\infty}^{x} \E_{\xi}T_{b}^{n-1} m(\xi)d\xi+\int_x^b (S(b)-S(\xi))\E_{\xi}T_{b}^{n-1} m(\xi)d\xi\right ).
\end{equation*}
If the integral in \eqref{eq:integral} diverges, using Fatou's lemma, we  have
$$\E_xT_b^n=n\int_{-\infty}^{+\infty}G(-\infty,b,x,\xi) \E_{\xi}T_{b}^{n-1} m(\xi)d\xi=\infty.$$
Hence independently of convergence or divergence of the integral \eqref{eq:integral} we have the equality \eqref{Kacb}.
The proof of \eqref{Kaca} is similar to this of \eqref{Kacb}. This finishes the proof of \eqref{Kacb} and \eqref{Kaca} for $n \in \N, $
$n\geq 1.$ 

We now turn to the proof of \eqref{Kacb} and \eqref{Kaca} for $p > 1, $ $ p \notin \N .$ Write $\alpha = p - [p]  \in ]0, 1 [ .$ Note that under our conditions, $[a, b] \ni x \mapsto \E_x T_{a,b}^\alpha $ is continuous which will be shown in lemma \ref{lem:cont} below.  
Hence exactly the same schema applies : We start from the function $f(x)=  \E_x T_{a,b}^{\alpha },$
  using \eqref{mainKac} we can write
  $$\E_xT^{1+\alpha}_{a,b}=(1+\alpha)\int_{-\infty}^{+\infty}G(a,b,x,\xi)\E_{\xi}T_{a,b}^{\alpha}m(\xi)d\xi.$$
  The function $x\mapsto \E_xT^{1+\alpha}_{a,b}$ is continuous on $[a,b],$ so we can apply the formula \eqref{mainKac} again. In each step we obtain a continuous function. Hence we can apply
  \eqref{mainKac}  $[p]$ times. In this way we obtain
  $$\E_xT^{p}_{a,b}=p\int_{-\infty}^{+\infty}G(a,b,x,\xi)\E_{\xi}T_{a,b}^{p-1}m(\xi)d\xi.$$
Then we pass to the limit when $a\to -\infty$ using exactly  the same considerations  as for \eqref {intKac}.
\end{proof}

The above proposition works for non-integer moments only if $[a, b] \ni x \mapsto \E_x T_{a,b}^\alpha $ is continuous, for any $0 < \alpha < 1.$ This is true under our conditions as shows the following lemma.

\begin{lem}\label{lem:cont}
Grant the assumption \ref{model}. Then $[a, b] \ni x \mapsto \E_x T_{a,b}^\alpha $ is continuous for any $ 0 < \alpha < 1 .$
\end{lem}

\begin{proof}
Let $x_n \to x , $ $x_n , x \in [a,b] .$ Write $F_n (dt)$
% = {\cal L} (T_{a,b} | P_{x_n} ), $
(respectively $ F(dt)$) for the law of $T_{a,b} $ under $\P_{x_n} $ (under $\P_x,$ respectively). Moreover, write 
$$ \varphi_n ( \lambda ) = \E_{x_n} e^{- \lambda T_{a,b} } , 
\;  \varphi  ( \lambda ) = \E_{x} e^{- \lambda T_{a,b} }$$
for the associated Laplace transforms.   
\begin{enumerate}
\item
We start by showing that for any $\lambda > 0 ,$ $ \varphi_n ( \lambda ) \to 
\varphi ( \lambda ) $ as $n \to \infty .$ 
For that sake, let $u_\lambda(x) $ for $a \le x \le b $ be the solution of the equation $$\left\{ \begin{array}{l}
L u_\lambda = \lambda u_\lambda \mbox{ in } ]a,b[ \\
u_\lambda(a) = u_\lambda (b) = 1 , 
\end{array}
\right\} .$$ 
Under our assumptions, the coefficients of the diffusion are H\"older-continuous on $[a,b] .$  
By continuity, $L$ is uniformly elliptic on $[a,b] .$ Hence, 
a solution to this problem exists and is given by
$$ u_\lambda(x) = \E_x e^{- \lambda T_{a,b} } ,$$
see \cite{KS}, chapter 5.7, proposition 7.2 and Remark 7.5,. This solution $u_\lambda(x)$ is continuous on $[a,b]$ which implies our claim. 
\item
By Maruyama and Tanaka, \cite{MaTa}, formula (3.7), we have that for any $0 < \alpha < 1$ and any $n,$
$$ \int_0^\infty \frac{1 - \varphi_n (\lambda)}{\lambda^{1 + \alpha}} d \lambda = \int_0^\infty t^\alpha F_n (dt) \int_0^\infty  \frac{1 - e^{-\lambda} }{\lambda^{1 + \alpha}} d \lambda .$$
In other words,
$$  \int_0^\infty \frac{1 - \varphi_n (\lambda)}{\lambda^{1 + \alpha}} d \lambda = \left( \E_{x_n} (T_{a,b}^\alpha) \right) \cdot \left( \int_0^\infty  \frac{1 - e^{-\lambda} }{\lambda^{1 + \alpha}} d \lambda \right) . $$
On the left hand side of the above formula we use dominated convergence.
Note that 
$$\frac{1 - \varphi_n (\lambda)}{\lambda} 
\le \E_{x_n} T_{a,b}  \le \sup_n \E_{x_n} T_{a,b} ,$$
where $ \sup_n \E_{x_n} T_{a,b}  $ is finite due to continuity on $[a,b]$ of the function $x\to\E_{x} T_{a,b},$ see 
(\ref{kac1}). Hence we can use the upper bound   
$$ \frac{1 - \varphi_n (\lambda)}{\lambda^{1 + \alpha}} 
\le \lambda^{- (1+ \alpha )} 1_{ [1, + \infty [ } ( \lambda) 
+  \sup_n \E_{x_n} (T_{a,b}) \lambda^{- \alpha} 1_{ [0,1] } ( \lambda).$$
Then by dominated convergence,
$$ \int_0^\infty \frac{1 - \varphi_n (\lambda)}{\lambda^{1 + \alpha}} d \lambda 
\to \int_0^\infty \frac{1 - \varphi  (\lambda)}{\lambda^{1 + \alpha}} d \lambda,$$
which in turn equals 
$$  \int_0^\infty \frac{1 - \varphi  (\lambda)}{\lambda^{1 + \alpha}} d \lambda = \left( \E_{x} (T_{a,b}^\alpha) \right) \cdot \left( \int_0^\infty  \frac{1 - e^{-\lambda} }{\lambda^{1 + \alpha}} d \lambda \right) , $$
applying once more formula (3.7) of \cite{MaTa}. This implies that 
$$  \E_{x_n} (T_{a,b}^\alpha) \to  \E_{x} (T_{a,b}^\alpha) \mbox{ as } n \to \infty ,$$
and this finishes our proof.
\end{enumerate}
\end{proof}

%Due to proposition \ref{prop:kac}, it is reasonable to introduce the following definition. 
%\begin{defin}\label{degreofrec}
%We say that the diffusion $X$ is positive recurrent of degree of recurrence  $p_0
%$, if
%$$p_0=\sup \{ p \geq 1  \quad \mbox{ such that for all }  x , y\in \R ,\quad \E_xT^p_y<\infty \}.$$ 
%\end{defin} 
%This definition has to be compared with the notion of positive recurrence, where for all $x,y$ $\E_xT_y<\infty$.
It is known, see for example \cite {MaTa}, that for $p>0$,  $x<b$ (resp $x>a$) the hitting time's moments satisfy the following property:  $\E_xT_b^p$ (resp $\E_xT_a^p$) is finite or infinite simultaneously for all couples $(x,b)$ s.t. $x<b$ (resp $(x,a)$ s.t. $x>a$ ). In the following theorem we refine this result and give an independent proof based on the generalized Kac's formula.

\begin{theo}\label{indepofpoint} 
Grant assumptions \ref{model} and \ref{posrec}. 
\begin{enumerate}
\item Let $x<b$ and $p \geq 1$.
\begin{description}
\item[(1i)] $\E_xT_b^p<\infty$ if and only if  $\int_{-\infty}^x \E_{\xi}T_b^{p-1}m(\xi)d\xi<\infty$.
\item [(1ii)] If for one couple $x<b,$ $\E_xT_b^p<\infty$, then for all couples $x'<b', $  $\E_{x'} T_{b'}^p<\infty.$ Moreover, for all
$b'$ fixed, the function $x'\mapsto \E_{x'} T_{b'}^p$ is continuous in $x',$ for $x'<b'$.
\end{description}

\item Let $a<x$ and $p \geq 1 $.
\begin{description}
\item[(2i)] $\E_xT_a^p<\infty$ if and only if  $\int_{x}^{+\infty} \E_{\xi}T_a^{p-1}m(\xi)d\xi<\infty$.
\item [(2ii)] If for one couple $a<x,$  $\E_xT_a^p<\infty$, then for all couples $a'<x', $  $\E_{x'} T_{a'}^p<\infty.$ Moreover, for all
$a'$ fixed, the function $x'\mapsto \E_{x'} T_{a'}^p$ is continuous in $x',$ for $a'<x'$.
%\item Let $n\in\N$ and $x,x',a,a'$ be four points in $E$ such that $a<x$ and $a'<x'$.
%If $\E_x T_a^n<\infty,$ then $\E_{x'} T_{a'}^n<\infty.$ 
\end{description}
\end{enumerate}
\end{theo}

\begin{proof}
1. Suppose $p=1.$ Using Kac's formula,
\begin{equation}\label {Kac1}
\E_xT_b
=(S(b)-S(x))\int_{-\infty}^xm(\xi)d\xi+\int_x^b (S(b)-S(\xi))m(\xi)d\xi . 
\end{equation}
The functions $S$ and $m$ are continuous, hence the last expression is finite if and only if $\int_{-\infty}^xm(\xi)d\xi<\infty$. The finiteness of the last integral does not depend on $x$ nor on $b .$ Hence, $\E_xT_b$  is finite or not simultaneously for all $x, b$ such that $x<b$. If  $\E_xT_b<\infty,$ the Kac's formula \eqref{Kac1} gives the continuity in $x<b$ of $\E_xT_b$ .

2. Now let $p = \alpha + 1 , $ where $\alpha \in ]0, 1[ . $ Suppose for some fixed $x<b, $ $\E_xT_b^p<\infty.$ Then $\E_x T_b < \infty,$ too. 
Hence $\E_{x'} T_{b'} < \infty $ for all $x' < b' .$ Then also $\E_{x'} T^\alpha_{b'} < \infty $ for all $x' < b' .$ By Kac's formula,
$$ \E_x T_b^{p} = p \left( (S(b)-S(x))\int_{-\infty}^x \E_\xi (T_b^\alpha ) m(\xi)d\xi+\int_x^b (S(b)-S(\xi)) \E_\xi (T_b^\alpha) m(\xi)d\xi \right) . $$
$\E_xT_{b}^p$ is finite if and only if $\int_{-\infty}^{x} \E_{\xi}T_{b}^{\alpha} m(\xi)d\xi<\infty.$
We have the upper bound $\E_\xi T_b^\alpha \le 1 + \E_\xi T_b ,$ where $ \xi \mapsto \E_\xi T_b$ has already been shown
to be continuous. Hence we see that for fixed $b$ the  integral $\int_{-\infty}^{x'} \E_{\xi}T_{b}^{\alpha } m(\xi)d\xi$ converges or diverges simultaneously  for all $x'<b.$ Hence we obtain the following equivalence for fixed $b\in \R .$
\begin{equation*}%\label{nequiv1}
\mbox {For some $x$ s.t. $ x<b$}  \; \E_xT_{b}^p<\infty\quad  \iff \quad\mbox { for all $ x'$ s.t. $  x'<b$}   \; \E_{x'}T_{b}^p<\infty . 
\end{equation*}
Then the continuity of $x \mapsto \E_x T_b^p ,$ $x < b,$ follows by dominated convergence, if $\E_x T_b^p < \infty.$  

Now let $\E_xT_{b}^p<\infty$ and fix some $b'$ such that $x<b<b'$.
We have $\E_xT_{b'}^p<\infty$ if and only if $\int_{-\infty}^{x} \E_{\xi}T_{b'}^{\alpha } m(\xi)d\xi<\infty.$
Using the strong Markov property and the sub-additivity of the function $x \mapsto x^{\alpha } ,$ we have 
 \begin{eqnarray*}%\label {markovestmation}
 \int_{-\infty}^{x} \E_{\xi}T_{b'}^{\alpha} m(\xi)d\xi &\leq &
 \left ( \int_{-\infty}^{x} \E_{\xi}T_{b}^{\alpha } m(\xi)d\xi+\E_{b}T_{b'}^{\alpha}\int_{-\infty}^{x}  m(\xi)d\xi\right )\\
 & \le &  \left ( \int_{-\infty}^{x} \E_{\xi}T_{b}^{\alpha } m(\xi)d\xi+[\E_{b}T_{b'} +1]\int_{-\infty}^{x}  m(\xi)d\xi\right ) .
 \end{eqnarray*}
Moreover, for $x<b<b'$, $\E_xT_{b}^p \leq \, \E_xT_{b'}^p.$
Therefore, the following two statements are equivalent. 
 \begin{equation*}
\mbox {For some $ x$ s.t. $x<b$}  \; \E_xT_{b}^p<\infty\quad \iff\quad\mbox {for all  $ b'$ s.t.$  x<b<b' $} \quad\E_{x}T_{b'}^p<\infty .
\end{equation*}
3. Now let $p \geq 2.$ We suppose the claim of the theorem verified for all moments of order $\alpha + k ,$ $1 \le k <   [p],$ and we show it for $p .$

Suppose for some fixed $x<b, $ $\E_xT_b^p<\infty.$ Then  $\E_xT_b^{p-1}<\infty,$ too. This implies
by our recurrence assumption that $\E_{x'}T_{b'}^{p-1}$ is finite and continuous for all $x'<b'.$  We use generalized Kac's formula once more in order to get
\begin{equation*}%\label{Kacn}
\E_xT_{b}^p=p\left( (S(b)-S(x))\int_{-\infty}^{x} \E_{\xi}T_{b}^{p-1} m(\xi)d\xi+\int_x^b (S(b)-S(\xi))\E_{\xi}T_{b}^{p-1} m(\xi)d\xi\right ).
\end{equation*}
 $\E_xT_{b}^p$ is finite if and only if $\int_{-\infty}^{x} \E_{\xi}T_{b}^{p-1} m(\xi)d\xi<\infty.$
Using continuity of $\E_{\xi}T_{b}^{p-1},$ we see that for fixed $b$ the  integral $\int_{-\infty}^{x'} \E_{\xi}T_{b}^{p-1} m(\xi)d\xi$ converges or diverges simultaneously  for all $x'<b.$
Hence we obtain the following equivalence for fixed $b\in \R .$
\begin{equation}\label{nequiv}
\mbox {For some $x$ s.t. $ x<b$}  \; \E_xT_{b}^p<\infty\quad  \iff \quad\mbox { for all $ x'$ s.t. $  x'<b$}   \; \E_{x'}T_{b}^p<\infty . 
\end{equation}
Now let $\E_xT_{b}^p<\infty$ and fix some $b'$ such that $x<b<b'$.
We have $\E_xT_{b'}^p<\infty$ if and only if $\int_{-\infty}^{x} \E_{\xi}T_{b'}^{p-1} m(\xi)d\xi<\infty.$
Using the strong Markov property and H\"older's inequality, 
 \begin{equation*}%\label {markovestmation}
 \int_{-\infty}^{x} \E_{\xi}T_{b'}^{p-1} m(\xi)d\xi\leq 2^{p-2} \left ( \int_{-\infty}^{x} \E_{\xi}T_{b}^{p-1} m(\xi)d\xi+\E_{b}T_{b'}^{p-1}\int_{-\infty}^{x}  m(\xi)d\xi\right ) .
 \end{equation*}
Moreover, for $x<b<b'$, $\E_xT_{b}^p \leq \, \E_xT_{b'}^p.$
Therefore, the following two statements are equivalent. 
 \begin{equation}\label{nnequiv}
\mbox {For some $ x$ s.t. $x<b$}  \; \E_xT_{b}^p<\infty\quad \iff\quad\mbox {for all  $ b'$ s.t. $  x<b<b' $} \quad\E_{x}T_{b'}^p<\infty .
\end{equation}
4. The proof of point $2.$ of the theorem is similar.
With \eqref{nequiv} and \eqref{nnequiv}, the proof is complete.
\end{proof}
\section{Estimation of moments for hitting times}
The question of existence of moments of hitting times arises in various problems and is widely  studied in the literature, see Fitzsimmons and Pitman \cite {FitsKac}, Carmona  and Klein \cite {CaKl}, Darling and Siegert \cite {DaSi}, Veretennikov \cite{Veret}, Balaji and Ramasubramanian \cite {bhatt}, Ditlevsen \cite {Di}, Deaconu and Wantz \cite {DeaWa}, Giorno, Nobile, Riccardi, Sacredote \cite {GNRS}, Kavian, Kerkyacharian, Roynette \cite {KKR} and the references therein. In this section we explore some sufficient and necessary conditions for existence of polynomial moments of hitting times and give lower and upper bounds on these moments. 

To give examples of diffusion with finite or infinite moments of hitting times, we have to impose some conditions on $\beta(x)$ and $\sigma^2(x)$ for large $|x|$. The first one guarantees the finiteness of the moments up to some order.
\begin{ass}\label{Hasm1}
There exist $M_0>0$, $\sigma_0>0$, $-\infty<\gamma<1$ and $r>0$ such that
\begin{equation*}%\label{sigma}
\sigma_0 |x|^\gamma\leq |\sigma(x)|\quad\mbox{and}\quad -\frac{x\beta(x)}{\sigma^2(x)}\geq r \quad\mbox {for}\ |x|>M_0. 
\end{equation*}
\end{ass}

\begin{ex}
This condition is for example satisfied for $\sigma(x)=1$ and $\beta(x) = -\frac{ \vartheta x}{1 + x^2},$ where $\vartheta > 1/2 ,$ for any $1/2<r < \vartheta$. For the recurrent Ornstein-Uhlenbeck process having $\beta (x) = - \vartheta x ,$ $\vartheta > 0 , $ Assumption \ref{Hasm1} is satisfied for any $r>1/2$ by taking $M_0$ large enough.
\end{ex}
It is well-known, see for instance \cite{bhatt}, that $\E_x T_a^n$ is finite for $n<r+1/2$ (if $\gamma=0$). However, in order to verify the conditions of our Theorems \ref{thm: Npolbound}, \ref{mainbounded} and \ref{mainint}, we have to estimate $\E_\nu T_a^{n}$ for $ n = p/2,$ thus the finiteness of $\E_x T_a^n$ is not sufficient for our purpose, we need a finer control on $\E_x T_a^n$ in order to control integrability of $\E_x T_a^n$ with respect to $\nu .$

The second Assumption, which is somewhat complementary to Assumption \ref{Hasm1}, ensures that starting from some order, the moments of hitting times are infinite. 
\begin{ass}\label{Hasm2}
There exist $M_0>0$, $\sigma_1>0$, $-\infty<\delta<1$ and $R>0$ such that
\begin{equation*}%\label{sigma}
0<|\sigma(x)|\leq\sigma_1 |x|^\delta \quad\mbox{and}\quad0<-\frac{x\beta(x)}{\sigma^2(x)}\leq R \quad\mbox {for}\ |x|>M_0. 
\end{equation*}
\end{ass}
%It is well-known, see also \cite {Veret}, that under condition \ref{Hasm2}, $\E_x T_a^n$ is infinite for $n>R/\sigma_0^2+1/2.$

\begin{ex}
We continue the previous example. For $\beta(x) = -\frac{ \vartheta x}{1 + x^2},$ where $\vartheta > 0 ,$ Assumption \ref{Hasm2} is satisfied for any $R > \vartheta .$ For the recurrent Ornstein-Uhlenbeck process, $\beta (x) = - \vartheta x,$ Assumption \ref{Hasm2} obviously does not hold.
\end{ex}

Note that the Assumptions \ref{Hasm1} and \ref{Hasm2} do not need to hold simultaneously.

Recall that the scale function of $X_t$ is given by
\[S(x)= \int_0^x s(t)dt,\quad \mbox{where}  \quad s(t)= \exp\left(-2\int_0^t\frac{\beta(u)}{\sigma^2(u)}du\right).\]
It is easily seen that under condition \ref{Hasm1}, $\int_{-\infty}s(x)\,dx=+\infty=\int^{\infty}s(x)\,dx$, which implies that $X_t$ is recurrent (not necessarily positive recurrent) with speed density $m(\xi)=\frac{1}{s(\xi)\sigma^2(\xi)}$. The speed density is precisely the density of the unique (up to a constant factor) invariant measure of the process, and positive recurrence is equivalent to the finiteness of the speed measure of~$X$, (see \cite{RY} Ch.VII, ¤3 and ¤4 ex .3.20, 3.21; \cite{BS} Ch. II ¤12).  %Note that under condition \ref{Hasm1}, for any $|x| \geq M_0,$ $$ s^{-1} (x) \le C \,  \left(\frac{M_0}{|x|}\right)^{2r} , \; m(x) \le C \frac{1}{\sigma_0^2 } \left(\frac{M_0}{|x|}\right)^{2r} ,$$ which is integrable if $ 2r > 1 .$ Hence under condition \ref{Hasm1}, the process is positive recurrent if $ r > 1/2 .$
In the sequel we will need to estimate the moments of the speed measure of $X$. We start with an elementary lemma.
\begin{lem}\label{lem:control}
Let $0<a\leq x$. Denote
\[I_{p,q}(x,a)=\int_x^\infty\frac{(\xi-a)^p}{\xi^q}\,d\xi\quad\mbox{and}\quad J_{p,q}(x,a)=\int_a^x\frac{(\xi-a)^p}{\xi^q}\,d\xi\]
Then
\[\frac{(x-a)^{p+1}}{(q-p-1)x^q}\leq I_{p,q}\leq \frac{x^{p+1}}{(q-p-1)x^q}\quad\mbox{for}\ 0\leq p<q-1\]
and
\[\frac{(x-a)^{p+1}}{\kappa x^q}\leq J_{p,q}\leq \frac{x^{p+1}}{(p+1-q)x^q}\quad\mbox{for}\ p\geq 0,\ q<p+1,\]
where $\kappa=(p+1)$ if $q>0$ and $\kappa=(p+1-q)$ if $q\leq 0$.
\end{lem}
\begin{proof}
Note that $I_{p,q}(x,a)<\infty$ if $0\leq p<q-1$. We have also $I_{p+1,q+1}(x,a)\leq I_{p,q}(x,a)$, whence
\begin{equation*}
I_{p,q}(x,a)=\left.\frac{(\xi-a)^{p+1}}{(p+1)\xi^q}\right|_x^\infty+\frac{q}{p+1}I_{p+1,q+1}(x,a)\leq\\
-\frac{(x-a)^{p+1}}{(p+1)x^q}+\frac{q}{p+1}I_{p,q}(x,a),
\end{equation*}
which yields
\[I_{p,q}(x,a)\geq \frac{(x-a)^{p+1}}{(q-p-1)x^q}.\]

In the same manner, for $p\geq 0$ and $q<p+1$,
\begin{equation*}
J_{p,q}(x,a)=\frac{(x-a)^{p+1}}{(p+1)x^q}+\frac{q}{p+1}J_{p+1,q+1}(x,a),
\end{equation*}
whence
\[J_{p,q}(x,a)\geq\left\{\begin{array}{rcl}
\frac{(x-a)^{p+1}}{(p+1)x^q}&\mbox{if}&q>0\\
\frac{(x-a)^{p+1}}{(p+1-q)x^q}&\mbox{if}&q\leq 0\ .
\end{array}\right.\]

On the other hand, under the respective conditions,
\[I_{p,q}(x,a)\leq \int_x^\infty\frac{\xi^{p}}{\xi^q}\,d\xi=\frac{x^{p+1}}{(q-p-1)x^q}\]
and
\[J_{p,q}(x,a)\leq \int_a^x\frac{\xi^{p}}{\xi^q}\,d\xi\leq\frac{x^{p+1}}{(p+1-q)x^q}\ .\]
\end{proof}

Put
\[p^*=\sup\{p>0\ :\ \int_x^{\infty}\frac{\xi^p}{\sigma^2(\xi)s(\xi)} d \xi<\infty\}.\]
Note that, for $M_0\leq x\leq\xi$,
\[\frac{s(x)}{s(\xi)}=\exp\left(2\int_x^\xi\frac{\beta(u)}{\sigma^2(u)}\right).\]
Under the assumption~\ref{Hasm1}, this yields $s(x)/s(\xi)\leq (x/\xi)^{2r}$ and $m(\xi)\leq C\xi^{-2r-2\gamma}$, hence $p^*\geq 2r+2\gamma-1$.

On the other hand, under the assumption~\ref{Hasm2}, $m(\xi)\geq C|\xi|^{-2R-2\delta}$, hence $p^*\leq 2R+2\delta-1$ and $\int_x^{\infty}\frac{\xi^{p^*}}{\sigma^2(\xi)s(\xi)} d \xi=\infty$ if $p^* = 2R+2\delta-1$.

\begin{theo}\label{theo:control}
Let $M_0<a<x$ or $x<a<-M_0$.
\begin{enumerate}
\item
Suppose that the assumption~\ref{Hasm1} holds with $2r+2\gamma>1$. For any positive real number $1\leq m<(2r+1)(1-\gamma)^{-1}/2$ put $\alpha=m-[m]$. Then 
\[\E_x T_a^m\leq \frac{x^{2m(1-\gamma)}}{r_m\sigma_0^{2m}(1-\gamma)^m}%+ o(x^{2m-1})
,\]
where $r_m=(2r+2\gamma-1)^\alpha\prod_{k=1}^{[m]} (2r-2(k+\alpha)(1-\gamma)+1)$.
%Hence $\E_x T_a^n<\infty$ for all $x$ and $a$ and $n<p^*/2+1$.

\item
Under the assumption~\ref{Hasm2}, for any integer $n\geq 1$:
\begin{itemize}
	\item if $n\leq p^*(1-\delta)^{-1}/2+1$ then
\[\E_x T_a^n\geq \frac{(x-a)^{2n(1-\delta)}}{R_n\sigma_1^{2n}\kappa^n}%+ o(x^{2m-1})
,\]
where $R_n=\prod_{k=1}^{n} (2R-2k(1-\delta)+1)$ and $\kappa=1\vee(1-\delta)$.
	\item if $n > p^*(1-\delta)^{-1}/2+1$, in particular if $n>(2R+1)(1-\delta)^{-1}/2$, then $\E_x T_a^n=\infty$.
\end{itemize}
\end{enumerate}
\end{theo}

\begin{rem}
Let us compare the above theorem to some known results. Note that most of them require that
\[0<\sigma_0^2\leq\sigma^2(x)\quad\mbox{or}\quad \sigma^2(x)\leq \sigma_1^2,\quad |x|>M_0,\]
i.e. $\gamma=\delta=0$ in our notations. To simplify the comparison, we assume it below, unless otherwise stated. Note, however, that our theorem holds under more general assumptions~\ref{Hasm1} and~\ref{Hasm2}.
\begin{enumerate}
\item
Under the condition 
\begin{equation*}%\label{sigma}
x\beta(x)<-r |x|^{1-p}  \quad\mbox {for}\ |x|>M_0 \mbox{ and } 0 < p < 1 ,
\end{equation*}
Douc, Fort and Guillin, \cite{DFG}, obtain the sub-exponential integrability of hitting times. They do not treat the critical case $p = 1$ which we consider here. 

\item It is known from Balaji and Ramasubramanian~\cite{bhatt} that, under the corresponding assumptions, $\E_x T_a^p<\infty$ for $p<r+1/2$ and $\E_x T_a^p=\infty$ for $p>R+1/2$. Nevertheless, they do note provide explicit bounds on $\E_x T_a^p$. Moreover, we show that in fact, at least for integer $n$, $\E_x T_a^n=\infty$ as soon as $n>p^*/2+1$, which can be much smaller then $R+1/2$.

\item
In~\cite[Theorem 4]{Veret}, Veretennikov obtains -- in our notations -- an upper bound 
\begin{equation}\label{veret1}
 \E_x T_a^p \le C ( 1 + |x|^m) \mbox{ for any $p < r_0 := r - 1 + \frac12  \frac{\sigma_0^2}{\sigma_1^2} $ }
\end{equation} 
and for any $m \in ] 2p, 2 r_0 [ .$ In our Theorem \ref{theo:control} we obtain an upper bound 
\begin{equation}\label{veret2}
\E_x T_a^p \le C( 1 + |x|^{2p} )\mbox{ for any $p < r + \frac12 .$ }  
\end{equation} 
Comparing \eqref{veret1} and \eqref{veret2}, we see that we have pushed the range of $p$ a little bit further: we obtain the control of moments up to at least $ p = r_0 + 1.$  Moreover, our constant $C$ for $x>a>M_0$ is quite explicit and sharp, as seen by taking a diffusion with constant drift and $r=R$.

%\item The lower bounds provided in Theorem \ref{theo:control} enable us to obtain necessary conditions for finiteness of $\E_\nu T_a^n$ for any measure $\nu$, and these results are new compared to the existing literature on hitting time moments.
\end{enumerate}
\end{rem}

\begin{proof}
Suppose $M_0\leq a<x$. If $\E_x T_a^m<\infty$ for $m\geq 1$, by the generalized Kac formula,
(theorem \ref{thm:mainKac})
\[\left(\E_x T_a^m\right)'=2m s(x)\int_x^{\infty}\frac{\E_\xi T_a^{m-1}}{\sigma^2(\xi)s(\xi)}\,d\xi , \]
where the derivative is taken with respect to $x$. Note that, for $x\leq\xi$,
\[\frac{s(x)}{s(\xi)}=\exp\left(2\int_x^\xi\frac{\beta(u)}{\sigma^2(u)}\right).\]

We start with the lower bound for $\E_x T_a^n$, $n\leq p^*(1-\delta)^{-1}/2+1$ under the assumption~\ref{Hasm2}. Note that the assertion is true if $\E_x T_a^n=\infty$, so we assume $\E_x T_a^n<\infty$ in the sequel. Recall that $p^*\leq 2R+2\delta-1$, so $n\leq (2R+1)(1-\delta)^{-1}/2$. Note also that
\[\frac{s(x)}{s(\xi)\sigma^2(\xi)}\geq \frac{x^{2R}}{\sigma_1^2\xi^{2R+2\delta}}.\]
Firstly suppose that $n<(2R+1)(1-\delta)^{-1}/2$. For $n=1$ we get, in the notations of lemma~\ref{lem:control}
\[\left(\E_x T_a\right)'=2s(x)\int_x^{\infty}\frac{d\xi}{\sigma^2(\xi)s(\xi)}\geq \frac{2x^{2R}}{\sigma_1^2}I_{0,2R+2\delta}(x,a)\geq \frac{2(x-a)}{(2R+2\delta-1)\sigma_1^2x^{2\delta}} , \]
whence
\[\E_x T_a\geq \frac{2}{(2R+2\delta-1)\sigma_1^2}\int_a^x \frac{\xi-a}{\xi^{2\delta}}\,d\xi=\frac{2}{(2R+2\delta-1)\sigma_1^2}J_{1,2\delta}(x,a)\geq
\frac{(x-a)^2}{R_1\sigma_1^2\kappa x^{2\delta}}\ .\]
By induction, for $n>1 ,$
\begin{multline*}
\left(\E_x T_a^n\right)'\geq \frac{2n s(x)}{R_{n-1}\sigma_1^{2n-2}\kappa^{n-1}}\int_a^x \frac{(\xi-a)^{2n-2}}{\xi^{2\delta(n-1)}\sigma^2(\xi)s(\xi)}\,d\xi\geq \frac{2n x^{2R}}{R_{n-1}\sigma_1^{2n}\kappa^{n-1}}I_{2n-2,2R+2\delta n}(x,a)\geq\\ \frac{2n(x-a)^{2n-1}}{R_{n-1}(2R+2\delta n-(2n-1))\sigma_1^{2n}\kappa^{n-1} x^{2\delta n}}= \frac{2n}{R_{n}\sigma_1^{2n}\kappa^{n-1}}\frac{(x-a)^{2n-1}}{x^{2\delta n}}, 
\end{multline*}
whence
\[\E_x T_a^n\geq \frac{2n}{R_{n}\sigma_1^{2n}\kappa^{n-1}}J_{2n-1,2\delta n}(x,a)\geq\frac{(x-a)^{2n}}{R_{n}\sigma_1^{2n}\kappa^{n}x^{2\delta n}} . \]

In the case $n=(2R+1)(1-\delta)^{-1}/2\leq p^*(1-\delta)^{-1}/2+1$, recalling that $p^*\leq 2R+2\delta-1$, we deduce that $p^*=2R+2\delta-1$ and $n=p^*(1-\delta)^{-1}/2+1$. Then 
\[\E_x T_a^{n-1}\geq C x^{2(n-1)(1-\delta)}=C x^{2R+2\delta-1}\]
for $x$ large enough. So $\int_x^\infty \E_\xi T_a^{n-1} m(\xi)\,d\xi=\infty$, and theorem~\ref{indepofpoint} yields that $\E_x T_a^{n}=\infty$. Hence the first point of the second assertion of the theorem is true for all $n\leq p^*(1-\delta)^{-1}/2+1$.

Now, suppose that $n\leq p^*(1-\delta)^{-1}/2+1 < n+1$, then $\E_x T_a^n\geq C x^{2n(1-\delta)}$ for $x$ large enough. Since $2n(1-\delta) > p^*$, we get $\int_x^\infty \E_\xi T_a^n m(\xi)\,d\xi=\infty$, and theorem~\ref{indepofpoint} yields again that $\E_x T_a^{n+1}=\infty$.

To prove the upper bound, note firstly that for any $p<2r+2\gamma-1$,
\[s(x)\int_x^{\infty}\frac{\xi^p\,d\xi}{\sigma^2(\xi)s(\xi)}\leq\frac{x^{2r}}{\sigma_0^2}\int_0^\infty \xi^{p-2r-2\gamma}\,d\xi=\frac{x^{p+1-2\gamma}}{\sigma_0^2(2r+2\gamma-p-1)}.\]
So we get
\[\left(\E_x T_a\right)'=2s(x)\int_x^{\infty}\frac{d\xi}{\sigma^2(\xi)s(\xi)}\leq \frac{2x^{1-2\gamma}}{(2r+2\gamma-1)\sigma_0^2},\]
whence
\[\E_x T_a\leq \int_a^x \frac{2\xi^{1-2\gamma}}{(2r+2\gamma-1)\sigma_0^2}\,d\xi=\frac{x^{2-2\gamma}-a^{2-2\gamma}}{r_1\sigma_0^2(1-\gamma)}\leq\frac{x^{2-2\gamma}}{r_1\sigma_0^2(1-\gamma)}.\]
Now, starting from
\[\E_x T_a^\alpha \le ( \E_x T_a)^\alpha \le \frac{x^{2\alpha(1-\gamma)}}{r_1^\alpha \sigma_0^{2\alpha}(1-\gamma)^\alpha} ,\]
we get analogously to the above calculus, applying successively the Kac formula :
\[ \E_x T_a^m \le \frac{x^{2m(1-\gamma)}}{r_m\sigma_0^{2m}(1-\gamma)^m}% + o ( x^{ 2m -1})
.\]

The case $x<a<-M_0$ follows by symmetry.
\end{proof}

\begin{rem}
Theorem~\ref{indepofpoint} implies the finiteness (and the continuity in $x$) or the infiniteness of $\E_x T^m_a$ for all $x$ and $a$ under the corresponding hypotheses of theorem~\ref{theo:control}.
\end{rem}

We would like to end this article with two corollaries, giving some ``practical'' form of the deviations theorems proved in section 3.
\begin{cor}\label{expbounded}
Suppose that $X$ satisfies the assumption \ref{model} and that the assumption~\ref {Hasm1} holds with $2r+2\gamma>1$. Take some $1<p<(2r+1)(1-\gamma)^{-1}/2$ and let $f\in\L^1(\mu)$, with $\|f\|_\infty <\infty$. Then for any initial distribution $\nu$ such that $\int_{\R} |x|^{p(1-\gamma)}d{\nu(x)}<\infty$, for all $0< \ge <  \|f\|_\infty $ and $t\geq 1$, the following inequality holds:  
\begin{equation}\label{polbounded}
\P_{\nu}\left(\left|\frac 1t\int_0^tf(X_s)ds-\mu(f)\right|>\ge\right)\leq 
\left\{
\begin{array}{ll}
K(l,p,\nu,X)\frac 1{\ge^p}\|f\|_\infty ^p \, {t^{- p/2}} & \mbox{ if } p \geq 2 \\
K(l,p, \nu,X)\frac 1{\ge^p}\|f\|_\infty ^p \, {t^{-\frac{p-1}{2}}} & \mbox{ if }1< p < 2
\end{array}
\right\}. 
\end{equation}
Here $K(l, p,\nu,X)$ is a positive constant, different in the two cases, which does not depend on $f$, $t$, $\ge$. In particular,~\eqref{polbounded} holds under $\P_x$ for all $x\in\R$.
\end{cor}
\begin{proof}
Assumption~\ref {Hasm1} together with $2r+2\gamma>1$ implies the positive recurrence of $X$ and also that $(2r+1)(1-\gamma)^{-1}/2>1$. Let $1<p<(2r+1)(1-\gamma)^{-1}/2.$ Since 
\[\E_{\nu}(R_1)^{p/2}\leq  (2^{p/2-1}\vee 1)(\E_{\nu}T_b^{p/2}+\E_bT_a^{p/2})\quad\mbox{and}\quad \E_{\nu}(R_2-R_1)^{p}=\E_a R_1^p,\]
we can see that the hypotheses of the theorem \ref{mainbounded} are satisfied if for some $a<b$ it holds that $\E_aT_b^p<\infty,$ $\E_bT_a^p<\infty$ and $\E_{\nu}T_b^{p/2}=\int_{\R} \E_xT_b^{p/2}\nu(dx)<\infty$. Using the theorem  \ref{theo:control} and the remark above, we obtain $\E_x T_y^p<\infty$ and $\E_x T_y^{p/2}<\infty$ for all $(x,y)\in\R^2$.%, which yields~\eqref{polbounded} under $\P_x$.

Without loss of generality we can choose $a=-M_0$ and $b=M_0.$ We then have
\[\E_x T_b^{p/2}\leq C{|x|^{p(1-\gamma)}},\ x>b\quad\mbox{and}\quad \E_x T_a^{p/2}\leq C{|x|^{p(1-\gamma)}},\ x<a.\]
Further, for any $x<a$, $\E_x T_b^{p/2}\leq C(\E_x T_a^{p/2}+\E_a T_b^{p/2})$. Finally, the continuity of $\E_x T_b^{p/2}$ implies
\[\E_x T_b^{p/2}\leq C{|x|^{p(1-\gamma)}}+C_1,\quad\mbox{for all}\ x\in\R,\] 
hence~\eqref{polbounded} holds under $\P_\nu$ if $\int_{\R} |x|^{p(1-\gamma)}d{\nu(x)}<\infty$.
\end{proof}

In the similar way one shows the following result:
\begin{cor}\label{explintegrable}
Suppose that $X$ satisfies the assumption \ref{model} and that the assumption~\ref {Hasm1} holds with $2r+4\gamma> 3$. Let $f $ be a bounded function with compact support. Let $\nu$ be an initial distribution such that $\int_{\R} |x|^{p(1-\gamma)}d{\nu(x)}<\infty$. Then for all $ p\in\N $, $2\leq p<(2r+1)(1-\gamma)^{-1}/2$, for all $0<\ge<\mu(|f|)$ and $t \geq 1$, the following inequality holds:
\begin{equation}\label{explpolintegr}
\P_{\nu}\left(\left|\frac 1t\int_0^tf(X_s)ds-\mu(f)\right|>\ge\right)\leq K(l,p,X )\frac 1{\ge^p}\mu(|f|)^p \, {t^{-p/2}}.
\end{equation}
Here $K(l, p , X)$ is a positive constant which does not depend on $f$, $t$, $\ge.$
\end{cor}

\def\refname{References}

\end{document}